\documentclass{article}

\usepackage[utf8]{inputenc}
\usepackage[margin=1in]{geometry}
\usepackage{graphicx}
\usepackage{mathrsfs}
\usepackage{wrapfig}
\usepackage{amsmath}
\usepackage{hyperref}
\usepackage{pgfgantt}
\usepackage{soul}
\usepackage{enumitem}
\usepackage{verbatim}
\usepackage{pdfpages}
\setlist[enumerate]{itemsep=2mm}

\usepackage{xcolor,colortbl}
\renewcommand*{\O}[1]{\mathcal{O}(#1)}

\newcommand*{\cfl}{cfl}

\newcommand*{\qxi}{\frac{\partial Q_i^n}{\partial x}}
\newcommand*{\qxxi}{\frac{\partial^2 Q_i^n}{\partial x^2}}
\newcommand*{\qxxxi}{\frac{\partial^3 Q_i^n}{\partial x^3}}
\newcommand*{\qti}{\frac{\partial Q_i^n}{\partial t}}
\newcommand*{\qtti}{\frac{\partial^2 Q_i^n}{\partial t^2}}
\newcommand*{\qttti}{\frac{\partial^3 Q_i^n}{\partial t^3}}
\newcommand*{\dx}{\Delta x}
\newcommand*{\dt}{\Delta t}

\graphicspath{{FIGURES/}}

\usepackage{authblk}
\title{Analysis of an Explicit, High-Order Semi-Lagrangian Nodal Method}
\author[1]{Gustaaf B. Jacobs\thanks{gjacobs@sdsu.edu}}
\author[1]{ Hareshram Natarajan}
\author[1]{ Pavel Popov}
\author[1,2]{ David. A. Kopriva }
\affil[1]{
Department of Aerospace Engineering and Computational Science Research Center, San Diego State University, San Diego, CA}
\affil[2]{Department of Mathematics,  Florida State University, Tallehassee, FL}
\date{}
\begin{document}

\maketitle

% $$
% \int_\Omega \left[ \sum_{i=0}^{N^{samples}} \frac{\partial u_i}{\partial t}+ \frac{W_i \partial u_i}{\partial x} = 0 =  \frac{\partial u}{\partial t}+ \frac{   \partial^2 D u}{\partial x^2} \right]
% $$

\begin{abstract}

A discrete analysis of the phase and dissipation errors of an explicit, semi-Lagrangian spectral element  method is performed.
The semi-Lagrangian method advects the Lagrange interpolant according the Lagrangian
form of the transport equations and uses a least-square fit to correct the update for interface 
constraints of neighbouring elements.
By assuming a monomial representation instead of the Lagrange form, a discrete version of the algorithm on a single element is derived. The resulting algebraic system lends itself to both a Modified Equation analysis and an eigenvalue analysis.
The Modified Equation analysis, which Taylor expands the stencil at a single
space location and time instance, shows that the semi-Lagrangian
method is consistent with the PDE form of the transport equation in the limit that the element size goes to zero.
The leading order truncation term of the Modified Equation is of the order of the degree of the interpolant which is consistent with numerical tests reported in the literature.
The dispersion relations 
show that the method is negligibly dispersive, as is common for semi-Lagrangian methods. 
An eigenvalue analysis shows that the semi-Lagrangian method with a nodal Chebyshev interpolant   
 is stable for a Courant-Friedrichs-Lewy condition based on the minimum collocation node spacing within an element that is greater than unity. \end{abstract}
% REQUIRED
%\begin{keywords}
%  Semi-Lagrangian, High-Order, Consistency, Accuracy, Stability, Modified Equation
%%\end{keywords}

% REQUIRED
%\begin{AMS}
%  65G99
%\end{AMS}
%
\section{Introduction}

In semi-Lagrangian (SL) methods, linear and non-linear transport equations of the form
\begin{equation}
    \frac{\partial q}{\partial t} + \frac{\partial f(q)}{\partial x} =0, 
    \label{eq:adv_pde}
\end{equation}
where the flux is $f(q)= a(q) q$ and $a(q)$ is the advection velocity,
are solved using the  Lagrangian form of the equivalent system of characteristic ordinary differential equations,
\begin{eqnarray}
    \frac{dx}{dt}&=&a (q), \label{eq:adv_x} \\
    \frac{dq}{dt} &=&0, 
    \label{eq:adv_ode}
\end{eqnarray}
which traces the solution along characteristics with a wave speed, $a(q)$, which is divergence
free, i.e., 
$\partial a(q)/\partial x= 0$ in one space dimension.
The discrete location at which solutions are approximated are traced according to  (\ref{eq:adv_ode}) over a single time step.
This is followed
by an interpolation of the advected solution back to its original location. This tracing
and mapping can be applied in a number of ways, using combinations of approximations of
(\ref{eq:adv_pde}), ((\ref{eq:adv_x}) and (\ref{eq:adv_ode}) and interpolation and numerical integration
methods  \cite{atm1,atm4,remap1,dgsl4,dgsl2}. 
%\gbj{cite semi-Lagrangian paper here}

In \cite{NatarajanJacobs20}, we developed an explicit, high-order,
semi-Lagrangian method that is consistent with a discontinuous spectral element method (DSEM) discretization
of a complementary Eulerian model that supplies the advection velocity $a(q)$. The semi-Lagrangian method approximates the solution on $P+1$ nodes supporting a degree of polynomial $P$ per element. The location and solution on each node is updated according to (\ref{eq:adv_ode}) using the Eulerian solver's $a(q)$. Combined with interface constraints that are found by connecting the solution on the interfaces of neighbouring elements, a least-square fit projects the solution back onto the original grid, where it can be coupled
back if desired to the Eulerian field solver.

 This semi-Lagrangian method ensures a boundary-fitted, higher-order accuracy approximation
 without the complications of Lagrangian-Eulerian, particle-mesh method interpolations
\cite{JH06,SJD14}. For further details on the method and its benefits, we refer the interested reader to \cite{NatarajanJacobs20}.
The algorithm is even more effective in the Monte-Carlo solution of Fokker-Planck type equations
using its equivalent of the Lagrangian-Langevin form \cite{NatarajanJacobs20}. A Lagrangian, Monte-Carlo update bypasses the curse-of-dimensionality
of the Fokker-Planck solution and the associated computational expense.

To connect the elemental solutions thus far, we have used an upwinding approach and interface constraints combined with the least-squares fit \cite{NatarajanJacobs20,NatarajanJacobs21}. 
In many DSEM solvers \cite{Kopriva09,HW08}  a Lax-Friedrichs (LF) method or other approximate Riemann solver is used to find the interface flux as a function of two neighbouring elemental solutions.
To be consistent with such upwinding, the LF coupling is thus desirable for the semi-Lagrangian update, too. An additional advantage of the LF solver
as compared to the upwind solver is that it is slightly easier to generalize for non-linear problems.

The impact of the interface flux specification on the numerical  properties
of the explicit, semi-Lagrangian method is not trivial. Semi-Lagrangian methods are known  for their
complete absence of dispersion as waves are traced exactly. Moreover, the method is stable for high
Courant-Friedrichs-Lewy ($\cfl$) numbers. This comes at a cost of considerable complexity of the traced basis and its mapping to the original basis.
In the explicit semi-Lagrangian method that we proposed in \cite{NatarajanJacobs20}, the characteristic information 
is only approximately enforced. As compared
to most semi-Lagrangian methods, this reduces the stable $cfl$ and may affect the diffusion and dispersion characteristics.

In addition to the basic properties of a numerical method, including accuracy, convergence, and stability, the properties of advection problem approximations are often also characterized through the diffusion and dispersion relations
that follow from an eigenvalue analysis of spatial flux derivative. An eigenvalue analysis provides a spectrum of the equivalent numerical wave frequencies as compared to the exact frequency, and shows to what extent waves at each frequency are accurately represented \cite{GassnerKopriva11}. 

Alternatively, for finite difference methods at least,
a Modified Equation Analysis (MEA) proves useful.
In MEA, the dependencies of the solution on multiple grid points within an  approximation stencil are reduced to a single nodal point using Taylor series, with the objective being to find an equivalent transport
PDE with a modified right hand side that summarizes the effect of truncation terms of the Taylor series. For higher-order methods, this analysis is not common. Only Costa \textit{et al.} \cite{CostaSherwinPeiro15}  report on MEA for modal discontinuous Galerkin methods. They present examples that provide evidence of the validity of a MEA analysis and its good comparison, albeit for lower wavenumbers, with an eigenvalue analysis \cite{GassnerKopriva11}.

Here, 
we revisit the semi-Lagrangian method of \cite{NatarajanJacobs20} and present a discrete analysis, including MEA  and an eigenvalue analysis of the recursive discrete system.
By assuming a monomial approximation in each element, we derive the discrete version of the algorithm on a single element with upwinding and Lax-Friedrichs interface patching. The resulting algebraic systems lend themselves to both the Modified
Equation Analysis (MEA) and an eigenvalue analyses.

The  approach we follow here is distinctly different from \cite{CostaSherwinPeiro15} for a number of reasons: firstly, we are concerned with a semi-Lagrangian method that does not require differentiation; secondly, we consider a collocation, strong form spectral method and thirdly, we use a monomial polynomial form for analysis. To the best of our knowledge, no such analysis has been attempted for high-order semi-Lagrangian methods.
The analysis sheds light on the influence of polynomial order, time stepping, grid spacing, interface patching  on dispersion, diffusion and stability of the method.
 
To introduce notation and for reference, we will first briefly review the explicit semi-Lagrangian method. Then we start with the derivation of the discrete zeroth order semi-Lagrangian method using the monomial form and analysis thereof. This sets stage for an extension to a first-order derivation and analysis followed
by a presentation for a general order $P$. 

\section{Semi-Lagrangian method}
\begin{figure}
    \centering
    \includegraphics[width=0.4\textwidth]{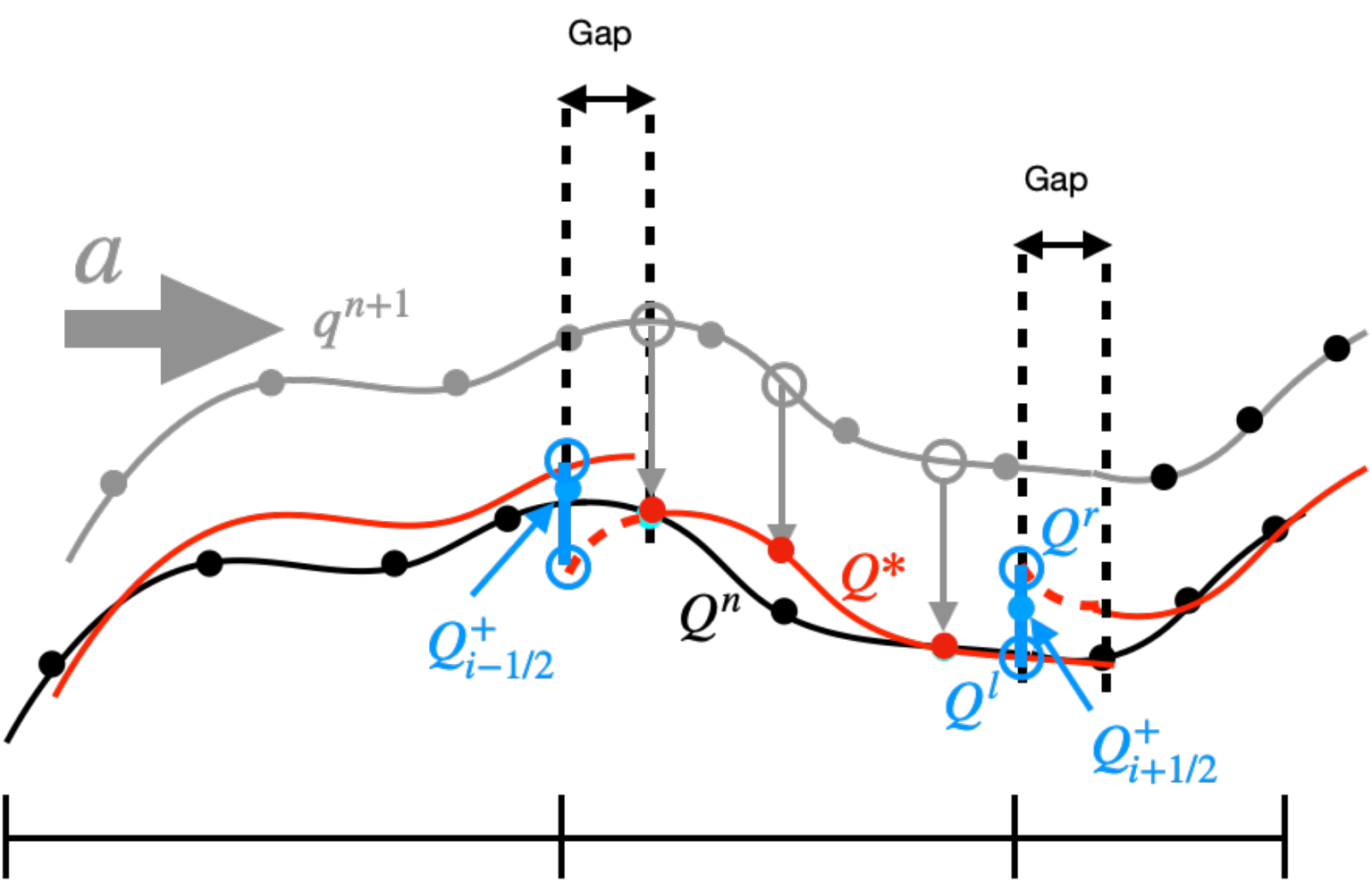}
    \caption{Schematic of the explicit, high-order Lagrangian method presented in \cite{NatarajanJacobs20}.}
    \label{fig:SLschematic}
\end{figure}

The  explicit, high-order semi-Lagrangian method  \cite{NatarajanJacobs20} solves the hyperbolic conservation laws in (\ref{eq:adv_pde})  on the domain $\Omega$. (See a schematic in Fig.\ \ref{fig:SLschematic}.) The method subdivides $\Omega$ into $K$ non-overlapping elements so that $\Omega=\sum_{k=1}^K \Omega_k$. Each element is mapped onto a local coordinate $\xi(x) \in [0,1]$. The element's solution at a given time $t^n$ is approximated with a Lagrange interpolant of $P^{th}$ degree as follows
\begin{eqnarray}
q_k^n(\xi) \approx Q^n(\xi) = \sum_{m=0}^P Q^n_m h_m(\xi),
\label{eq:SLN_initial_lagr}
\end{eqnarray}
where $Q_m^n$ are the solution values at nodal points $\xi_m$ with superscript $n$ identifying the time instance. 
The basis function, $h_m(\xi)$, is the Lagrangian interpolating polynomial defined at the Chebyshev Gauss quadrature nodes $\xi_m$, 
\begin{equation} %\label{eq:gaussquad}
\xi_{m} = \frac{1}{2}\left[ 1- \cos\left(\frac{m+1/2}{P}\right) \pi \right] \qquad m=0,1,...,P
\label{eq:nodalpoints}
\end{equation}
mapped onto $[0,1]$.

The nodal coordinates and solution values, i.e., the interpolant, are advected in time according to (\ref{eq:adv_x}) and (\ref{eq:adv_ode}) 
\begin{eqnarray}
 \xi_m^\dagger & =&  \xi_m^n + \Delta t a,  \nonumber\\
 Q^\dagger_{m} & = & Q^n_{m}
\end{eqnarray}
with an explicit time integrator (Euler for the sake of simplicity of notation here).
The advected interpolant is then projected back to the quadrature nodes so that 
\begin{eqnarray}
 Q^\star(\xi) =  \sum_{m=0}^P Q^\dagger_m h_m^\dagger(\xi)=\sum_{m=0}^P Q^\star_m h_m(\xi),
\label{eq:SLQstar}
\end{eqnarray}
with $h_m^\dagger(\xi)$ the Lagrange polynomial defined on $\xi_m^\dagger$.

It remains to correct the solution for neighbouring elements and boundary conditions. To identify boundary locations, we introduce a global counter, $i$, whose coordinate, $x_i$, coincides with the center coordinate of element $k$, i.e., $\xi(x_i)=1/2$. Variables with fractional subscript counters $i+1/2$ identify inter-element face locations so that 
$\Omega_k\in[x_{i-1/2}, x_{i+1/2}]$.  For the MEA below, this notation is particularly helpful, as we will find a Modified Equation at the global coordinate $x_i$ using Taylor expansions.
The boundary value at the interface between element $k$ and $k+1$ is denoted as $Q_{i+1/2}^+$, see Fig. \ref{fig:SLschematic}, which is determined from the face values of neighbouring elements, which we denote as $Q^l$  (to the left of the interface) and $Q^r$ (to the right of the interface). We therefore write
\begin{equation}
    Q_{i+1/2}^+ = \mathcal{G}
    \left(Q^l, Q^r\right),
%\left(q^*_{k}(x_{1+1/2}), q^*_{k+1}(x_{1+1/2}\right), 
    \label{eq:riemann_sol}
\end{equation}
 where the superscript $+$ identifies the inter-element values that results from  the solution of a Riemann problem according to some function $\mathcal{G}$ from the two face values. 
 $Q^l$ and $Q^r$ are interpolated to the edge locations $x_{i+1/2}$ using (\ref{eq:SLQstar}).
 
  In our previous work \cite{NatarajanJacobs20} we took $\mathcal{G}$ to be the upwind solution.
 Alternatively, here we will consider a more general Lax-Friedrichs weighting of the solutions from either side \cite{leveque_2002}
\begin{equation}
    \mathcal{G}
% \left(q^*_{k}(x_{1+1/2}), q^*_{k+1}(x_{1+1/2}\right)) = \frac{q^*_{k}(x_{1+1/2}) +q^*_{k+1}(x_{1+1/2})}{2} + \frac{\Delta t\|a\|_\infty}{2 \Delta x} \left( f(q^*_{k}(x_{1+1/2}) - f(q^*_{k+1}(x_{1+1/2}) \right),
\left(Q^l, Q^r\right) = \frac{Q^l +Q^r}{2} + \frac{\Delta t}{2 \Delta x} \left( f(Q^l) - f(Q^r) \right),
    \label{eq:LF_method}
\end{equation}
where 
\begin{equation}
\Delta t = \cfl \frac{\|\Delta x\|_{min}}{\|a\|_\infty}, 
\label{eq:cfl}
\end{equation}
and $\|a\|_\infty$ is the maximum wave speed in $\Omega$ and $\|\Delta x\|_{min}$ the minimum
spacing between the nodal points in the global coordinate $x$. The Courant-Lewy-Friedrichs number is denoted by $\cfl$.

With the two additional constrains on each face of the element, a  least-squares minimization is used to update the solution to time level $t^{n+1}$
as
\begin{eqnarray}
Q^{n+1}(\xi) = \sum_{m=0}^P Q^{n+1}_m h_m(\xi) = \arg \min_{f(\xi) \in \Re^P}  \left\|   \mathbf{Q}^{b} - f(\xi) \right\|_2,
\end{eqnarray}
where $\mathbf{Q}^{b}$ is the vector that contains
the nodal solution values  $\mathbf{Q}^\star = [ Q^\star_0,  Q^\star_1, \ldots  \quad, Q^\star_{P}]$ 
and the boundary constraints,
\begin{equation}
    \mathbf{Q}^{b} =\left[ Q_{i-1/2}^+, \mathbf{Q}^\star, Q_{i+1/2}^+ \right].
\label{eq:SLL_Qb}
\end{equation}

\section{Analysis of the Semi-Lagrangian Approximation}

In this section, we derive the modified equation associated with the semi-Lagrangian approximations. We start with the lowest order, $P=0$, for its simplicity, and then increase the complexity by analysing $P=1$ before moving to the general case. The results for the low order schemes give insight into the more general results, and the derivations are easier to follow.

\subsection{Zeroth order, $P$=0}

First, we analyze the discrete  semi-Lagrangian method for a  zeroth order ($P=0$) 
 approximation of the solution within the $k^{th}$ element at its center nodal point, $i$,
\begin{equation}
   Q_i   = q_k(x_i). 
\end{equation}
The nodal coordinate,  $x_i$, and solution are first updated according to (\ref{eq:adv_x}) and (\ref{eq:adv_ode}),
\begin{eqnarray}
    x_i^{n+1} &=& x_i^n+ \Delta t a^{n+1} \\
    q_k(x_i^{n+1})  &=& Q_i^n,
\end{eqnarray}
where the time step is chosen so that the nodal point remains within the element's bounds 
$x_{i-1/2} < x_i^{n+1} < x_{i+1/2}$,
\begin{equation}
    \Delta t \le cfl \Delta x/a,
    \label{eq:cfl_P0}
\end{equation}
with $cfl \le 0.5$ and $\Delta x = x_{i+1/2}-x_{i-1/2}$.

The solution is then interpolated back to the nodal point $x_i$.
Because the approximation is zeroth order and because the solution remains within the bounds of the element,
the updated, interpolated, solution at the nodal point is the same as the solution at $t_n$,
\begin{equation}
    Q_i^\star = Q_i^n.
    \label{eq:SL_projected}
\end{equation}
The interpolant is corrected for boundary constraints, $Q_{i-1/2}^\star$  and $Q_{i+1/2}^\star$ at both edges of the element using a least
square approach, so that the error
\begin{equation}
    \epsilon(Q_i^{n+1}) = \left( Q_{i-1/2}^\dagger -Q_i^{n+1}\right)^2 + \left( Q_{i}^\star -Q_i^{n+1}\right)^2 + \left( Q_{i+1/2}^\dagger -Q_i^{n+1}\right)^2
\end{equation}
is minimized. 
Setting the derivative $\frac{d  \epsilon(Q_i^{n+1})}{d Q_i^{n+1} }=0$, it follows
that 
\begin{equation}
    Q_i^{n+1} = ( Q_{i-1/2}^\dagger +  Q_{i}^\star  + Q_{i+1/2}^\dagger )/3.
    \label{eq:SL_LS}
\end{equation}
% \gbj{could do analysis on this using first order approximation of $Q_i^{n+1}  = ax +b$ and minimize with respect to $(a,b)$. 
% More algebra, but very doable and maybe worthwhile for us to understand the SL method better.}

It remains to determine the boundary constraints. To connect neighbouring volumes,
the Riemann problem
\begin{equation}
    Q_{i+1/2}^\star = \mathcal{G}\left(Q_{i}, Q_{i+1}\right)
    \label{eq:SL0_riemann_sol}
\end{equation}
 must be solved or approximated.
We propose to use a Lax-Friedrichs upwinding as
\begin{equation}
  Q_{i+1/2}^\dagger = \frac{Q_{i+1}^n + Q_{i}^n}{2} - \frac{3 \Delta t}{2 \Delta x} \left( f(Q^n_{i+1}) - f(Q^n_{i}) \right) ,
  \label{eq:SL0_bc}
\end{equation}
where, comparing to (\ref{eq:LF_method}), we have multiplied the flux term with a factor three
so that it cancels out against the denominator in (\ref{eq:SL_LS}).
Substituting \eqref{eq:SL0_bc} and (\ref{eq:SL_projected}) into (\ref{eq:SL_LS})  yields
\begin{eqnarray}
    Q_i^{n+1} & = &  \frac{1}{3}\left\{ \frac{Q_{i}^n + Q_{i-1}^n}{2} - \frac{3 \Delta t}{2 \Delta x} \left( f(Q^n_{i}) - f(Q^n_{i-1}) \right) \right. \\
              & + & Q_i^n \nonumber \\
              & + & \left. \frac{Q_{i+1}^n + Q_{i}^n}{2} - \frac{3 \Delta t}{2 \Delta x} \left( f(Q^n_{i+1}) - f(Q^n_{i}) \right) \right\} \nonumber
\end{eqnarray}
Simplifying leads to
\begin{equation}
    Q_i^{n+1}  =    \frac{ Q_{i-1}^n+4 Q_i^n +   Q_{i+1}^n}{6} -  \frac{\Delta t}{2 \Delta x}  \left( f(Q^n_{i+1}) - f(Q^n_{i-1})  \right).
    \label{eq:SL0_final}
\end{equation}

\textit{Remarks:}

When $P = 0$,
\begin{itemize}
 \item The stencil in (\ref{eq:SL0_final}) is similar to the Lax-Friedrichs (LF) method in (\ref{eq:LF_method}), but the average term on the right-hand side is distinctly different. Like LF, the method is diffusive, which can be understood if 
    $\frac{1}{3} Q_i^n$ is subtracted from both sides of (\ref{eq:SL0_final}). Then a second order central difference approximation to a diffusive term emerges with a  diffusion coefficient of $(\Delta x ^2)/3$. For LF, a similar term can be shown to appear, but with a diffusion coefficient of $(\Delta x ^2)/2$, making LF more diffusive. 
     \item The factor three in (\ref{eq:SL0_bc}) is required to obtain the second term on the right hand-side in (\ref{eq:SL0_final}), which correctly represents the approximation to the flux derivative in the Eulerian form and thus, in light of the previous remark, ensures conservation and a correct wave speed of the approximate solution.  As opposed to Eulerian DSEM methods, the semi-Lagrangian method is not flux based and thus is not trivially conservative. Because of the flux correction inherent to the Lax-Friedrichs method, the flux difference in (\ref{eq:SL0_final}) recovers the conservation property. We will explore this flux correction for general high-order semi-Lagrangian schemes in the sections below.
    \item If a simple upwind method where used in (\ref{eq:SL0_riemann_sol}) instead of LF, then that would eliminate the flux term and wave speed from the approximation and thus result in a meaningless scheme.
    \item Because the solution in (\ref{eq:SL_projected}) does not change at the nodal location for  the zeroth order
    advection and projection,  this step can be considered Eulerian in some sense.
    \item A least-squares method belongs to a general class of minimization methods that are commonly used to approximate Eulerian
    forms of the equation using Galerkin methods.
    \end{itemize}
\subsubsection{The Modified Equation}
Assuming linear advection with constant $a$, and Taylor expanding all terms in (\ref{eq:SL0_final}) around $x_i$ and $t^n$,  we find the ME
    \begin{equation}
        \frac{\partial Q^n_{i}}{\partial t} +  \frac{\partial ( a Q^n_{i})}{\partial x} = \left(\frac{1}{3} \frac{(\Delta x)^2}{\Delta t} - \frac{1}{2} (\Delta t)  a^2 \right) \frac{\partial^2 Q^n_i}{\partial x^2}  + \left(\frac{1}{6} (\Delta x)^2 a + \frac{1}{6}(\Delta t)^2 a ^3 \right) \frac{\partial^3 Q^n_i}{\partial x^3}+...
    \label{eq:SL_LF_equivalent}
    \end{equation}
Here, we have used (\ref{eq:adv_pde}) to relate temporal derivatives to spatial derivatives. The ME shows that the leading order term is $\O{\Delta x}$ and implies that the $P$=0 semi-Lagrangian method is first order accurate.
    
The ME confirms that the approximation is diffusive and dispersive. For stability of the scheme, it is required that
the diffusion coefficient be positive, which can be ensured by a time step restriction
    \begin{equation}
        \Delta t \le \left(\sqrt{\frac{2}{3}} \frac{(\Delta x)}{a}\right) .
    \end{equation}
Computational tests (not presented here) of a linear advection of sine wave show that this stability criterion  is not sufficient. For $cfl>$0.6, numerical solutions are unstable. {A Von Neumann analysis on the discrete stencil (\ref{eq:SL0_final}) confirms this more precisely, and that the stable $cfl<\frac{1}{\sqrt{3}}\approx 0.577$. } Note, that  the time step $\Delta t$ should obey (\ref{eq:cfl}) with $cfl<0.5$ because the $x_i$ should not be advected beyond the right edge of an element.
    
    % Interesting results as it show that the semi-Lagrangian method with a flux correction has a time step restriction and that the restriction is less than the first-order Eulerian Lax-Friedrichs scheme (factor is  $\sqrt{2/3}$ instead of 1/2 for Euler). 
    
   Rewriting (\ref{eq:SL_LF_equivalent}) in terms of $cfl$ using (\ref{eq:cfl}) leads to
    \begin{equation}
        \frac{\partial Q_{i}}{\partial t} +  \frac{\partial ( a Q_{i})}{\partial x} = \frac{\Delta x a}{6cfl}\left(2-3cfl^2 \right) \frac{\partial^2 Q_i}{\partial x^2}  + \frac{1}{6} (\Delta x)^2 a (1+cfl^2)    \frac{\partial^3 Q_i}{\partial x^3}+...,
    \end{equation}
    which shows that both the dissipation and dispersion error will reduce with grid spacing. The $cfl$ in the denominator of the first term on the right hand side means that the numerical diffusion will go to infinity in the limit of a zero time step, while maintaining a constant grid spacing. To mitigate this, the time and grid spacing must be reduced simultaneously. It is likely that this is a property is inherent to the Eulerian-Lagrangian duality of a semi-Lagrangian method.

\subsection{First-Order, $P$=1, Semi-Lagrangian Method}

To derive the discrete form of the semi- Lagrangian method for $P>0$, we find that a monomial representation of the basis enables analysis that the Lagrangian interpolant in (\ref{eq:SLN_initial_lagr}) does not. To derive the discrete stencil, we therefore start from a first-order monomial defined on two nodal points, $\xi_1$ and $\xi_2$ within the $k^{th}$ element as
\begin{equation}
    q_k^n(x) \approx  Q^n(x) =   C_1^n + C_2^n x,
\end{equation}
where $Q^n(x)$ is the polynomial of degree one at $t^n$.
The coefficient $C_2$ is the slope of the line and is determined as
\begin{equation}
    C_2^n = \frac{Q_1^n-Q_0^n}{\xi_1-\xi_0}.
    \label{eq:SL_C1}
\end{equation}
It then follows that $C_1$ is
\begin{eqnarray}
  C_1^n = \frac{Q_0^n \xi_1  - Q_1^n \xi_0}{\xi_1-\xi_0}.
  \label{eq:SL1_C2}
\end{eqnarray}
\noindent Once again, in the first step of the semi-Lagrangian method, the polynomial is advected according to its advection velocity. For simplicity we take a constant $a$, so that 
\begin{eqnarray}
 {\xi_0}^{\dagger} &=& {\xi_0}^n + a\Delta t , \nonumber \\
 {\xi_1}^{\dagger} &=& {\xi_1}^n + a\Delta t,\\
 Q_0^{\dagger} &=& Q_0^n \nonumber, \\
  Q_1^{\dagger} &= &Q_1^n,
\end{eqnarray}
which translates the straight line in the positive $x$-direction. A linear interpolation projects this line back onto the
original grid as
\begin{equation}
  Q^\star(x)  = Q^n(x) - \left(\frac{Q_1^n - Q_0^n}{\xi_1 - \xi_0} \right) a\Delta t.
  \label{eq:SL1_projected}
\end{equation}
We can rewrite \eqref{eq:SL1_projected} using (\ref{eq:SL_C1}) as
\begin{equation}
  Q^\star(x)  = Q^n(x) - C_2^n a\Delta t.  
  \label{eq:SL1_projected2}
\end{equation}
Note that the second term on the right hand side  of (\ref{eq:SL1_projected2}) represents the solution shift of the linear approximation with an element from 
$t^n$ to $t^{n+1}$ and is thus per the Eulerian conservation law  indicative of the mass flux through the element interface over the time step $\Delta t$ .
% from $x_{i-1/2}$  to $x_{i+1/2}$ represents
% the change in the mass flux within the element over time $\Delta t$.
In fact, we can use the definition of the flux $F(Q)=aQ$ to rewrite the equation in terms of the flux difference as
\begin{equation}
  Q^\star(x)  = Q^n(x) - \Delta t \frac{F_1^n - F_0^n}{\xi_1 - \xi_0}.
  \label{eq:SL1_projecteda}
\end{equation}
\noindent This $Q^{\star}$, together with boundary constraints $Q^+_{i-1/2}$ and $Q^+_{i+1/2}$, determines
$Q^{n+1}(x)$ using a least-square fit. For this fit we define the error as:
\begin{eqnarray}
  \epsilon(C_2^{n+1}, C_1^{n+1}) &=& (Q_{i-1/2}^+ - C_2^{n+1} x_{i-1/2} - C_1^{n+1})^2 + \nonumber \\  
   && (Q_{0}^\star \ \ \ \ \ - C_2^{n+1} \xi_0 - C_1^{n+1})^2  +\nonumber  \\
    && (Q_{1}^\star \ \ \ \ \ - C_2^{n+1} \xi_1 - C_1^{n+1})^2 + \nonumber  \\
     && (Q_{i+1/2}^+ - C_2^{n+1} x_{i+1/2} - C_1^{n+1})^2. \nonumber  \\
\end{eqnarray}
\noindent Minimization with respect to $C_1^{n+1}$ and $C_2^{n+1}$ leads to 
\begin{eqnarray}
  - 2E_1  + E_3 C_2^{n+1}   + E_4 C_1^{n+1}  & = &  0 \nonumber \\
  - 2E_2  + E_4 C_2^{n+1}   + 8  C_1^{n+1}  & = &  0,
 \label{eq:SL1_LSstep1}
\end{eqnarray}
where 
\begin{eqnarray} 
  E_1  &=&  \left(Q_{i-1/2}^+ x_{i-1/2} + Q_{0}^\star  \xi_0+   Q_{1}^\star \xi_1 +  Q_{i+1/2}^+ x_{i+1/2} \right), \nonumber \\
  E_2  &=&  \left(Q_{i-1/2}^+  + Q_{0}^\star +   Q_{1}^\star +  Q_{i+1/2}^+ \right), \nonumber \\
  E_3  & = &  2 \left( {x_{i-1/2}}^2 + {\xi_0}^2 + {\xi_1}^2 + {x_{i+1/2}}^2 \right), \nonumber \\
  E_4  & = &  2 \left( {x_{i-1/2}} + {\xi_0} + {\xi_1} + {x_{i+1/2}} \right).
   \label{eq:SL1_LSstep2}
\end{eqnarray}
Solving for  $C_2^{n+1}$ and $C_1^{n+1}$ yields
\begin{eqnarray}
 C_2^{n+1} & = & \frac{2E_1   - E_4 C_1^{n+1}}{E_3} = \frac{2E_1(8E_3 -E_3 E_4) - 2E_2 E_3 E_4+2E_1 E_4}{8E_3^2-E_3^2E_4},\\
  C_1^{n+1} & = & \frac{2E_2E_3-2E_1E_4}{8E_3-E_3E_4}.
  \label{eq:SL1_C1C2}
\end{eqnarray}
Note that $E_3$ and $E_4$ are predetermined and depend only on the coordinates.
Only $E_1$ and $E_2$ depend on the solution. $Q_1^\star$ and $Q_2^\star$ have a flux update
according to (\ref{eq:SL1_projected}).  
% We will take $x_i$=0 and assume the nodal points are symmetric around $x_i$=0. 

\begin{figure}[h] 
\centering 
\mbox{ \includegraphics[width=0.5\textwidth]{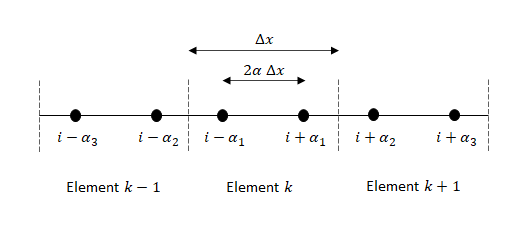}
      }
\caption{Schematic showing the 1D grid for the second order SL method with a spacing of $2 \alpha \dx$ between the nodes.} 
\label{fig:alpha_grid}
\end{figure}

 To understand the impact of the distribution of the nodal points within an element, we will consider
 a variable, symmetric nodal point distribution per element per the schematic shown in Fig. \ref{fig:alpha_grid}. We will assume the spacing between $\xi_1$ and $\xi_2$ to be $2 \alpha \Delta x$ so that  $\xi_1 = x_{i-\alpha_1}=-\alpha \Delta x$ and $\xi_2 = x_{i+\alpha_1}= \alpha \Delta x$. 
 Then we can write out $E_1$ and $E_2$ and simplify using (\ref{eq:SL1_projected}) to get the interior point values
\begin{eqnarray}
  Q_{i+\alpha_1}^\star &=& Q_{i+\alpha_1}^n - \Delta t \frac{F_{i+\alpha_1}^n-F_{i-\alpha_1}^n}{2 \alpha \Delta x} \nonumber \\
                       &=& Q^n_{i+\alpha_1} \left(1- a\frac{\Delta t }{2 \alpha \Delta x} \right) +
                           Q^n_{i-\alpha_1} \left(a\frac{\Delta t }{2 \alpha \Delta x} \right), \nonumber \\
  Q_{i-\alpha_1}^\star &=& Q_{i-\alpha_1}^n - \Delta t \frac{F_{i+\alpha_1}^n-F_{i-\alpha_1}^n}{2 \alpha \Delta x} \nonumber \\
                       &=&  Q^n_{i+\alpha_1} \left(-a\frac{\Delta t }{2 \alpha \Delta x} \right) + Q^n_{i-\alpha_1} \left(1+ a\frac{\Delta t }{2 \alpha \Delta x} \right).
  \label{eq:SL1_star}
  \end{eqnarray}
  For the symmetric node distribution, $E_4$=0 and $E_3=\Delta x^2(1/2 + 2 \alpha)$, so (\ref{eq:SL1_C1C2}) can be further simplified to
\begin{eqnarray}
 C_2^{n+1} & =& \frac{2 E_1}{\Delta x^2(1/2 + 2 \alpha)}, \\
 C_1^{n+1} &=&  \frac{E_2}{4}.
 \label{eq:SL1_C1C2_mod}
 \end{eqnarray}
 
With these coefficients the updated monomial at $t^{n+1}$ is determined. What remains are the boundary constraints in $E_1$ and $E_2$. We will discuss using an upwind scheme and a LF scheme for the boundary constraints in the next two subsections.

%%%%%%%%%%%%%%%%%%% Figure 1%%%%%%%%%%%%%%%%%%%%%%%%%%%%%

% General spacing
\subsubsection{The Modified Equation}

In this section we consider two interface conditions, one with pure upwinding, and the other with the LF coupling.

{\bf Upwinding.} Assuming a positive constant advection velocity, $a>0$, the upwind scheme generates the interface constraints, $Q_{i\pm 1/2}^+$, as 
% \begin{eqnarray}
% Q_{i-1/2}^+         &=& q^*_{k-1}(x_{1-1/2}), \\ 
% Q_{i+1/2}^+         &=& Q^*_{i+1/2}, 
% \end{eqnarray}
\begin{eqnarray}
Q_{i-1/2}^+         &=& Q^r_{k-1} = \quad Q^n_{i-\alpha_2} \left[\frac{1}{2} + \frac{1}{4\alpha} - a\frac{\Delta t }{2 \alpha \Delta x} \right] + Q^n_{i-\alpha_3} \left[\frac{1}{2} - \frac{1}{4\alpha} + a\frac{\Delta t }{2 \alpha \Delta x} \right], \nonumber \\ 
Q_{i+1/2}^+         &=& Q^r_k \quad =  \quad Q^n_{i+\alpha_1} \left[\frac{1}{2} + \frac{1}{4\alpha} - a\frac{\Delta t }{2 \alpha \Delta x} \right] + Q^n_{i-\alpha_1} \left[\frac{1}{2} - \frac{1}{4\alpha} + a\frac{\Delta t }{2 \alpha \Delta x} \right],
\label{eq:SL1_upwind_interface}
\end{eqnarray}
where the subscript $k$ identifies the element number. 
% Using Equation \ref{eq:SL1_star} and $2 \alpha \Delta x$ as the spacing between the points,
% \begin{eqnarray}
% Q_{i+\alpha_1}^\star &=& Q^n_{i+\alpha_1} \left(1- a\frac{\Delta t }{2 \alpha \Delta x} \right) +
%                         Q^n_{i-\alpha_1} \left(a\frac{\Delta t }{2 \alpha \Delta x} \right), \\ 
% Q_{i-\alpha_1}^\star &=& Q^n_{i+\alpha_1} \left(-a\frac{\Delta t }{2 \alpha \Delta x} \right) +
%                         Q^n_{i-\alpha_1} \left(1+ a\frac{\Delta t }{2 \alpha \Delta x} \right).
% \end{eqnarray}
Then using (\ref{eq:SL1_C1C2}) , the updated solution at the global grid coordinate, $i$, is given by
\begin{eqnarray}
Q_i^{n+1} &=& C_1^{n+1} + C_2^{n+1}(0), \nonumber \\
&=& \frac{1}{4}\left( Q_{i-1/2}^+ + Q_{i-\alpha_{1}}^\star + Q_{i+\alpha_{1}}^\star + Q_{i+1/2}^+ \right).
\label{eq:SL1_Q-i-n}
\end{eqnarray}
\noindent Substituting (\ref{eq:SL1_star}) and (\ref{eq:SL1_upwind_interface}) into \eqref{eq:SL1_Q-i-n} leads to the discrete stencil
\begin{eqnarray}
Q_i^{n+1} &=& Q^n_{i+\alpha_1} \left[\frac{3}{8} + \frac{1}{16 \alpha} - \frac{3a \Delta t}{8 \alpha \Delta x}         \right] +
Q^n_{i-\alpha1} \left[ \frac{3}{8} - \frac{1}{16 \alpha} + \frac{3a \Delta t}{8 \alpha \Delta x}         \right] + \nonumber \\
&&  Q^n_{i-\alpha_2} \left [\frac{1}{8} + \frac{1}{16 \alpha} - \frac{a \Delta t}{8 \alpha \Delta x}    \right] +
  Q^n_{i-\alpha_3} \left [\frac{1}{8} - \frac{1}{16 \alpha} + \frac{a \Delta t}{8 \alpha \Delta x}    \right].
\label{eq:SL1_stencil_upwind}
\end{eqnarray}
To find the Modified Equation, each term is Taylor expanded around $x_i$ and $t^n$,
\begin{eqnarray}
Q_i^{n+1}         &=& Q_i^n + \Delta t \qti + \frac{\Delta t^2}{2}\qtti + \frac{\Delta t^3}{6}\qttti + ... , \nonumber \\
Q_{i+\alpha_1}^{n} &=& Q_i^n + \alpha \Delta x \qxi +  \frac{1}{2}\alpha^2 \Delta x^2 \qxxi + \frac{1}{6}\alpha^3 \Delta x^3 \qxxxi + ... , \nonumber \\ 
Q_{i-\alpha_1}^{n} &=& Q_i^n - \alpha \Delta x \qxi +  \frac{1}{2}\alpha^2 \Delta x^2 \qxxi - \frac{1}{6}\alpha^3 \Delta x^3 \qxxxi + ... , \nonumber \\ 
Q_{i-\alpha_2}^{n} &=& Q_i^n - (1-\alpha) \Delta x \qxi +  \frac{1}{2}(1-\alpha)^2 \Delta x^2 \qxxi - \frac{1}{6}(1-\alpha)^3 \Delta x^3 \qxxxi + ... , \nonumber \\ 
Q_{i-\alpha_3}^{n} &=& Q_i^n - (1+\alpha) \Delta x \qxi +  \frac{1}{2}(1+\alpha)^2 \Delta x^2 \qxxi - \frac{1}{6}(1+\alpha)^3 \Delta x^3 \qxxxi + ... 
\label{eq:SL1_TaylorExpansions}
\end{eqnarray}   
From the PDF from of the advection equations (\ref{eq:adv_pde}), it follows that  $\qtti= a^2 \qxxi$ and $\qttti=- a^3 \qxxxi$.  Using these substitutions and substituting the Taylor expansions \eqref{eq:SL1_TaylorExpansions} into (\ref{eq:SL1_stencil_upwind}) yields the ME,
% \gbj{copy the result from the symbolic script here}
% \haresh{included the result}
\begin{eqnarray} 
 \qti +  a\qxi &=& \qxxi\left[\frac{1}{2}(\frac{\alpha^2\dx^2}{\dt}) +\frac{1}{4} (a\dx)-\frac{1}{2}(a^2\dt)   \right] + \nonumber \\
 & & \qxxxi \left[\frac{1}{48} (\frac{\dx^3}{\dt})  - \frac{1}{6}(a^3\dt^2)  - \frac{1}{12}(\frac{\alpha^2 \dx^3}{\dt}) - \frac{1}{8}(a \dx^2)  - \frac{1}{6}(a \alpha^2 \dx^2) \right]+...
 \label{eq:SL2_upwind_general} 
\end{eqnarray}
Replacing $\Delta t$ according to (\ref{eq:cfl}) with $\|\Delta x\|_{min}= \Delta x \left( \frac{1}{2} - \alpha \right)$, we can rewrite the modified equation as
\begin{eqnarray}
\qti +  a\qxi &=& a\dx\qxxi \left[ \frac{-2\alpha^2 -2 \cfl^2(\alpha -1/2)^2 - \cfl(\alpha-1/2)}{4\cfl(\alpha -1/2)}   \right] + \nonumber \\
               & & a\dx^2\qxxxi \left[ \frac{1}{24 \cfl(1-2 \alpha)}  - \frac{ \cfl^2 (1-2\alpha)^2}{24}  - \frac{\alpha^2}{6(1-2\alpha)} - \frac{  1}{8}  - \frac{ \alpha^2 }{6} \right]+ ...
\end{eqnarray}
% \qti +  a\qxi &=& a\dx\qxxi \left[ \frac{ \alpha^2 }{ \cfl (1-2\alpha)}   - \frac{1}{4}  \cfl (1- 2\alpha) + \frac{1}{4}   \right] + \nonumber \\
%               & & a\dx^2\qxxxi \left[ \frac{1}{24 \cfl(1-2 \alpha)}  - \frac{ \cfl^2 (1-2\alpha)^2}{24}  - \frac{\alpha^2}{6(1-2\alpha)} - \frac{  1}{8}  - \frac{ \alpha^2 }{6} \right]+ ...
The leading order term is of $\O{\Delta x}$, i.e., the $P=1$ method is also first order accurate. The diffusion and dispersion
coefficients both depend on the $\cfl$ and the grid point distribution, $\alpha$.
% The diffusion coefficient, $\nu$ is,
% \begin{eqnarray}
%         \nu = \frac{1}{2}\frac{\alpha^2\dx^2}{\dt} +\frac{1}{4} (a\dx)-\frac{1}{2}(a^2\dt).
% % \end{eqnarray}
% For the scheme to be stable the diffusion coefficient be greater than unity:
% \begin{eqnarray}
%          \frac{1}{2}\frac{\alpha^2\dx^2}{\dt} +\frac{1}{4} (a\dx)-\frac{1}{2}(a^2\dt) \geq 0
% \end{eqnarray}

% We can replace $\Delta t$ using (\ref{eq:cfl}) with
% If we use a $cfl$ number based on $\dx_{\text{min}}$ then the condition becomes, 
% \begin{eqnarray}
%     \dt &=& cfl \left( \frac{\dx_{\text{min}}}{a} \right), \nonumber \\
%     \dt &=& cfl \left( \frac{\dx}{a}(\frac{1}{2}-\alpha) \right).
% \end{eqnarray}

% -(a*(      ))/(48*cfl*dx*(alfa - 1/2))

{\bf Lax-Friedrichs.}
% \gbj{copy the result from the symbolic script here }
% \haresh{Added the result: $\omega$ should be cfl* $\omega$, may be use $\lambda$ for cfl number?}
% Similar to the previous section on the upwind fluxes, we consider the spacing between $\xi_1$ and $\xi_2$ to be $2 \alpha \Delta x$ and the points indexed as $\xi_1 = x_{i-\alpha_1}$ and $\xi_2 = x_{i+\alpha_1}$ respectively. The schematic of the grid is shown in Figure \ref{fig:alpha_grid}. 
% The extension of the SL solver with LF interface patching should be straightforward. As compared to the SL algorithm
% that was presented in \cite{} only the boundary constraints should be changed according to (\ref{eq:SL0_bc}).
% Note that in the DSEM approach the grid spacing and time step are not as distinctly defined
% as in the finite volume case. To circumvent the usage of these variable, we replace
% the time step in (\ref{eq:SL0_bc}) according to the stable  time step 
% \begin{equation}
%     \Delta t = cfl \Delta x/ a,
% \end{equation}
% where $cfl$ is the CFL condition, so that
% \begin{equation}
%   Q_{i+1/2}^\star = \frac{Q_{i+1}^n + Q_{i}^n}{2} -  \frac{ \omega cfl }{2 a} \left( f(Q^n_{i+1}) - f(Q^n_{i}) \right).
%   \label{eq:SL0_bc}
% \end{equation}
% We have removed the factor three as that is a factor is specific to the first-order method and replaced it with a factor $\omega$ which may vary per order of the scheme. We for now assume
% that $\omega$=1 for the general case.
The Lax-Friedrichs scheme determines the interface constraint according to (\ref{eq:LF_method}), as 
\begin{eqnarray}
Q_{i-1/2}^+         &=&  \frac{Q^r_{k-1} +Q^l_k}{2} + \frac{\omega \Delta t}{2 \Delta x} \left( f(Q^r_{k-1}) - f(Q^l_k) \right), \nonumber \\ 
Q_{i+1/2}^+         &=& \frac{Q^r_{k} +Q^l_{k+1}}{2} + \frac{\omega \Delta t}{2 \Delta x} \left( f(Q^r_{k}) - f(Q^l_{k+1}) \right),
\label{eq:SL1_LF}
\end{eqnarray}
where we have introduced a factor $\omega$  to adjust the weight of the flux correction, similar to the factor of three that was used in the zeroth order method above. 
%An increased magnitude of the factor is well to known to enhance diffusion.
Note that for a constant advection velocity and with $F=aQ$, the upwind method in (\ref{eq:SL1_upwind_interface}) is recovered if $\omega=\frac{\Delta x }{ a \Delta t } $, and thus the ME should be the same for LF as compared to upwind for this $\omega$.

Assuming a positive constant advection velocity, $a>0$, the LF scheme according to (\ref{eq:SL1_LF}) determines the interface constraints, $Q_{i\pm 1/2}^+$, as
% \gbj{haresh, to be consistent with the upwind section, you need to write out the equivalent discrete version of (\ref{eq:SL1_upwind_stencil}) and (\ref{eq:SL1_stencil_upwind}) for Lax-Friedrichs here. Please do so.}
% \haresh{I have added the discrete stencil terms for the LF method}
\begin{eqnarray}
  % Q_{i-1/2}^+ &=& \quad  \ \ \left( Q^n_{i-\alpha_3} \left[\frac{1}{2} - \frac{1}{4\alpha} + a\frac{\Delta t }{2 \alpha \Delta x} \right] +  Q^n_{i-\alpha_2} \left[\frac{1}{2} + \frac{1}{4\alpha} - a\frac{\Delta t }{2 \alpha \Delta x} \right] \right) \left(\frac{1}{2} + \frac{\omega a \dt}{2\dx} \right)+  \nonumber \\
  %       & & \quad  \ \ \left( Q^n_{i-\alpha_1} \left[\frac{1}{2} + \frac{1}{4\alpha} + a\frac{\Delta t }{2 \alpha \Delta x} \right] +  Q^n_{i+\alpha_1} \left[\frac{1}{2} - \frac{1}{4\alpha} - a\frac{\Delta t }{2 \alpha \Delta x} \right] \right) \left(\frac{1}{2} - \frac{\omega a \dt}{2\dx} \right) \nonumber     
  Q_{i-1/2}^+ &=& \quad  \ \ Q^n_{i-\alpha_3} \left[\frac{1}{4} - \frac{1}{8\alpha} + \frac{a \Delta t }{4 \alpha \Delta x} + \frac{\omega a \dt}{4 \dx} - \frac{\omega a \dt}{8 \alpha \dx} + \frac{\omega a^2 \dt^2 }{4 \alpha \dx^2}\right] + \nonumber \\ 
            & & \quad \ \ Q^n_{i-\alpha_2} \left[\frac{1}{4} + \frac{1}{8\alpha} - \frac{a \Delta t }{4 \alpha \Delta x} + \frac{\omega a \dt}{4 \dx} + \frac{\omega a \dt}{8 \alpha \dx} - \frac{\omega a^2 \dt^2 }{4 \alpha \dx^2} \right]+   \nonumber \\ 
            & & \quad  \ \ Q^n_{i-\alpha_1} \left[\frac{1}{4} + \frac{1}{8\alpha} + \frac{a \Delta t }{4 \alpha \Delta x} - \frac{\omega a \dt}{4 \dx} - \frac{\omega a \dt}{8 \alpha \dx} - \frac{\omega a^2 \dt^2 }{4 \alpha \dx^2} \right]+   \nonumber \\
            & & \quad  \ \ Q^n_{i+\alpha_1} \left[\frac{1}{4} - \frac{1}{8\alpha} - \frac{a \Delta t }{4 \alpha \Delta x} - \frac{\omega a \dt}{4 \dx} + \frac{\omega a \dt}{8 \alpha \dx} + \frac{\omega a^2 \dt^2 }{4 \alpha \dx^2} \right]   \nonumber \\
   Q_{i+1/2}^+ &=& \quad  \ \ Q^n_{i-\alpha_1} \left[\frac{1}{4} - \frac{1}{8\alpha} + \frac{a \Delta t }{4 \alpha \Delta x} + \frac{\omega a \dt}{4 \dx} - \frac{\omega a \dt}{8 \alpha \dx} + \frac{\omega a^2 \dt^2 }{4 \alpha \dx^2}\right] + \nonumber \\ 
            & & \quad \ \ Q^n_{i+\alpha_1} \left[\frac{1}{4} + \frac{1}{8\alpha} - \frac{a \Delta t }{4 \alpha \Delta x} + \frac{\omega a \dt}{4 \dx} + \frac{\omega a \dt}{8 \alpha \dx} - \frac{\omega a^2 \dt^2 }{4 \alpha \dx^2} \right]+   \nonumber \\ 
            & & \quad  \ \ Q^n_{i+\alpha_2} \left[\frac{1}{4} + \frac{1}{8\alpha} + \frac{a \Delta t }{4 \alpha \Delta x} - \frac{\omega a \dt}{4 \dx} - \frac{\omega a \dt}{8 \alpha \dx} - \frac{\omega a^2 \dt^2 }{4 \alpha \dx^2} \right]+   \nonumber \\
            & & \quad  \ \ Q^n_{i+\alpha_3} \left[\frac{1}{4} - \frac{1}{8\alpha} - \frac{a \Delta t }{4 \alpha \Delta x} - \frac{\omega a \dt}{4 \dx} + \frac{\omega a \dt}{8 \alpha \dx} + \frac{\omega a^2 \dt^2 }{4 \alpha \dx^2} \right].   
\label{eq:SL1_LF_stencil}
\end{eqnarray}

Substituting (\ref{eq:SL1_star}) and (\ref{eq:SL1_LF_stencil}) into \eqref{eq:SL1_Q-i-n} leads to the discrete stencil
\begin{eqnarray}
 Q_{i}^{n+1} &=& \quad \ \  Q^n_{i+\alpha_1} \left[\frac{3}{8}  \ \    \quad \quad  \quad - \frac{3    a\dt}{8\alpha \dx} \  + \omega \frac{a \Delta t  }{ \Delta x} \left( \quad \quad \ \  + \frac{ 1}{16 \alpha }  \quad \quad  \quad \quad \right)\right] + \nonumber \\
 &&  \quad \ \  Q^n_{i-\alpha_1} \left[ \frac{3}{8} \ \   \quad \quad\quad  + \frac{3 a\dt}{8\alpha \dx}    \ + \omega \frac{a \Delta t  }{ \Delta x} \left(  \quad \quad   \ \    - \frac{ 1}{16 \alpha}   \quad \quad \quad  \quad \right) \right]+ \nonumber \\
            & & \quad \ \ Q^n_{i-\alpha_2} \left[\frac{1}{16} + \frac{1}{32\alpha} - \frac{a \Delta t }{16 \alpha \Delta x}   + \omega  \frac{a \Delta t }{\Delta x}\left(  \ \  \frac{ 1}{16 } + \frac{ 1}{32 \alpha } - \frac{ a \dt }{16 \alpha \dx}\right)  \right]+   \nonumber \\
            & & \quad  \ \ Q^n_{i-\alpha_3} \left[   \frac{1}{16} - \frac{1}{32\alpha} + \frac{a \Delta t }{16 \alpha \Delta x}   + \omega   \frac{a \Delta t}{ \Delta x} \left(  \ \ \frac{ 1}{16 } - \frac{ 1}{32 \alpha} + \frac{ a \dt }{16 \alpha \dx}\right)\right] + \nonumber \\ 
            & & \quad  \ \ Q^n_{i+\alpha_2} \left[\frac{1}{16}  + \frac{1}{32\alpha}  + \frac{ a \Delta t }{16 \alpha \Delta x} + \omega   \frac{a \Delta t}{\Delta x} \left(-\frac{ 1}{16 } - \frac{ 1}{32 \alpha } - \frac{ a \dt }{16 \alpha \dx} \right)\right]+   \nonumber \\
            & & \quad  \ \ Q^n_{i+\alpha_3} \left[\frac{1}{16}  - \frac{1}{32\alpha}  - \frac{ a \Delta t }{16 \alpha \Delta x}  + \omega  \frac{a \Delta t}{\Delta x} \left( -\frac{ 1}{16} + \frac{ 1}{32 \alpha } + \frac{ a \dt }{16 \alpha \dx} \right)\right].
            \label{eq:SL2_LF_stencil}
 \end{eqnarray} 
 Note that if we substitute $\omega=\frac{\Delta x }{ a \Delta t } $,
 the upwind stencil in (\ref{eq:SL2_upwind_general}) is consistently recovered.

 \newcommand*{\qkpaa}{Q_{i+\alpha_1}^n}
 \newcommand*{\qkpbb}{Q_{i+\alpha_2}^n}
 \newcommand*{\qkpcc}{Q_{i+\alpha_3}^n}
  \newcommand*{\qkmaa}{Q_{i-\alpha_1}^n}
 \newcommand*{\qkmbb}{Q_{i-\alpha_2}^n}
 \newcommand*{\qkmcc}{Q_{i-\alpha_3}^n}

Using the Taylor's series expansion and the relations between space and time derivatives from the PDE, we arrive
at the ME for the  semi-Lagrangian scheme with the Lax-Friedrichs interface constraint,
\begin{eqnarray}
    \qti  + a\qxi &=& \qxxi\left[\frac{1}{2} (\frac{\alpha^2 \dx^2}{\dt}) + \frac{1}{4}(a^2 \omega \dt) - \frac{1}{2} (a^2 \dt)   \right]+  \nonumber \\
    & & \qxxxi \left[ \frac{1}{48}(a \dx^2 \omega) -\frac{1}{6} (a^3 \dt^2) - \frac{1}{12}(a \alpha^2 \dx^2 \omega) - \frac{1}{8} (a \dx^2)  - \frac{1}{6} (a \alpha^2 \dx^2)    \right]+... .
 \label{eq:SL2_LF_general} 
\end{eqnarray}
This scheme, like the scheme with the upwind constraint per (\ref{eq:SL2_upwind_general}), is therefore first order accurate according to the leading order term on the right-hand side. 
Consistently, if we compare the two MEs for the upwind and LF constraint, i.e.,   (\ref{eq:SL2_upwind_general}) and (\ref{eq:SL2_LF_general}), then
it is clear  that once again with $\omega=\frac{\Delta x}{ a \Delta t}$ the MEs are the same.
Replacing $\Delta t$ according to (\ref{eq:cfl}) with $\|\Delta x\|_{min}= \Delta x \left( \frac{1}{2} - \alpha \right)$, we can rewrite \eqref{eq:SL2_LF_general} as
\begin{eqnarray}
\qti  + a\qxi &=& a \dx \qxxi
\left[ \frac{- 2 \alpha^2 +  (2-\omega)\cfl^2  (\alpha - 1/2)^2}{4 cfl (\alpha -1/2)} \right] + \nonumber \\
 & &  a \dx^2 \qxxxi \left[ - \frac{cfl^2(\alpha - 1/2)^2}{6} - \frac{1}{8}   + \frac{\omega}{48} -  \frac{\alpha^2}{6} - \frac{\alpha^2 \omega}{12}\right] +... .
\label{eq:SL2_LF_general_replace_dt} 
\end{eqnarray}

We explore the stability and accuracy characteristics of \eqref{eq:SL2_LF_general_replace_dt} by setting
the diffusion coefficient to zero, which requires that
\begin{eqnarray}
   2 - \omega - \frac{2\alpha^2}{cfl^2(\alpha^2- 1/2)}  =0.
    \label{eq:SL2_zerodiffusion}
\end{eqnarray}
If the diffusion coefficient is greater than zero the method
is diffusive in the second order term. For a zero diffusion coefficient, the method increases in accuracy by one order and the diffusion is in the fourth order derivative term with respect to $Q$. Just like for the zeroth order method discussed in the previous section, it is possible
to set $\omega$ so that the accuracy increases, but for $P=1$ this depends on $\alpha$, i.e.,
the distribution of the nodal locations, as well. Thus, for any given $\alpha$ and $cfl$, there is a corresponding
$\omega$ that sets the diffusion coefficient to zero. 

\begin{figure}[h] 
\centering 
\mbox{ \includegraphics[width=0.45\textwidth]{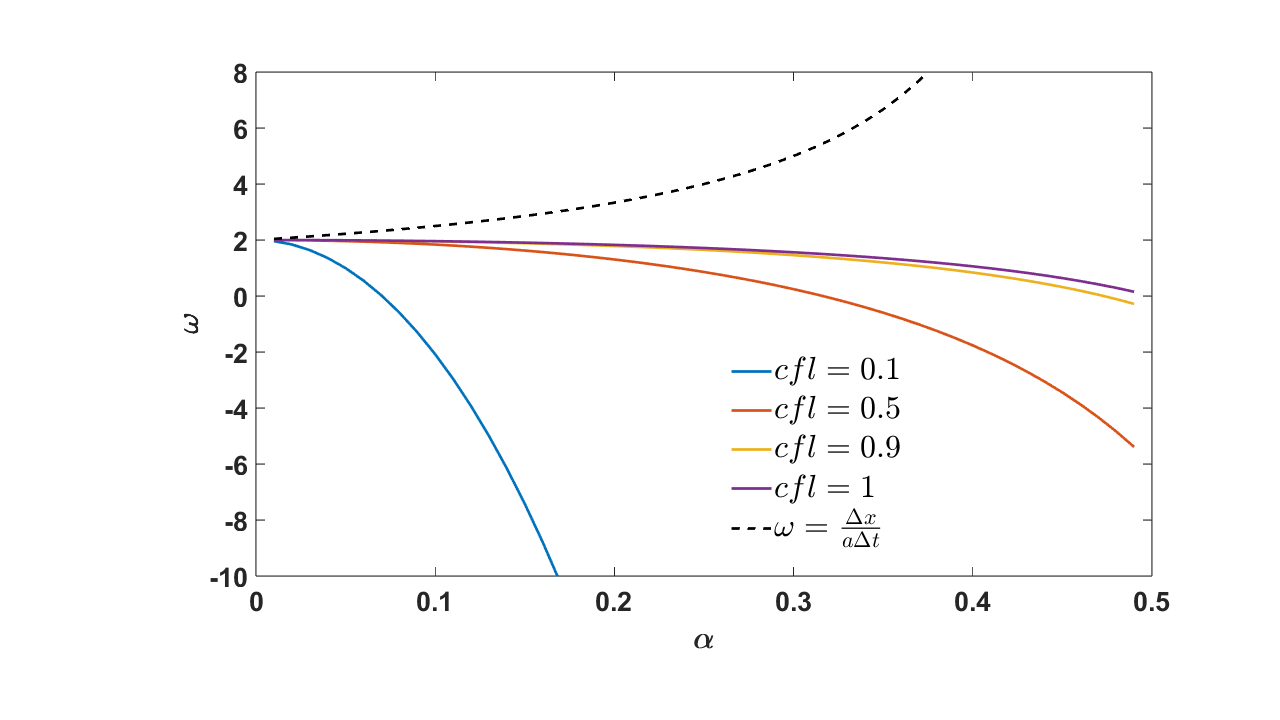}      }
\mbox{ \includegraphics[width=0.45\textwidth]{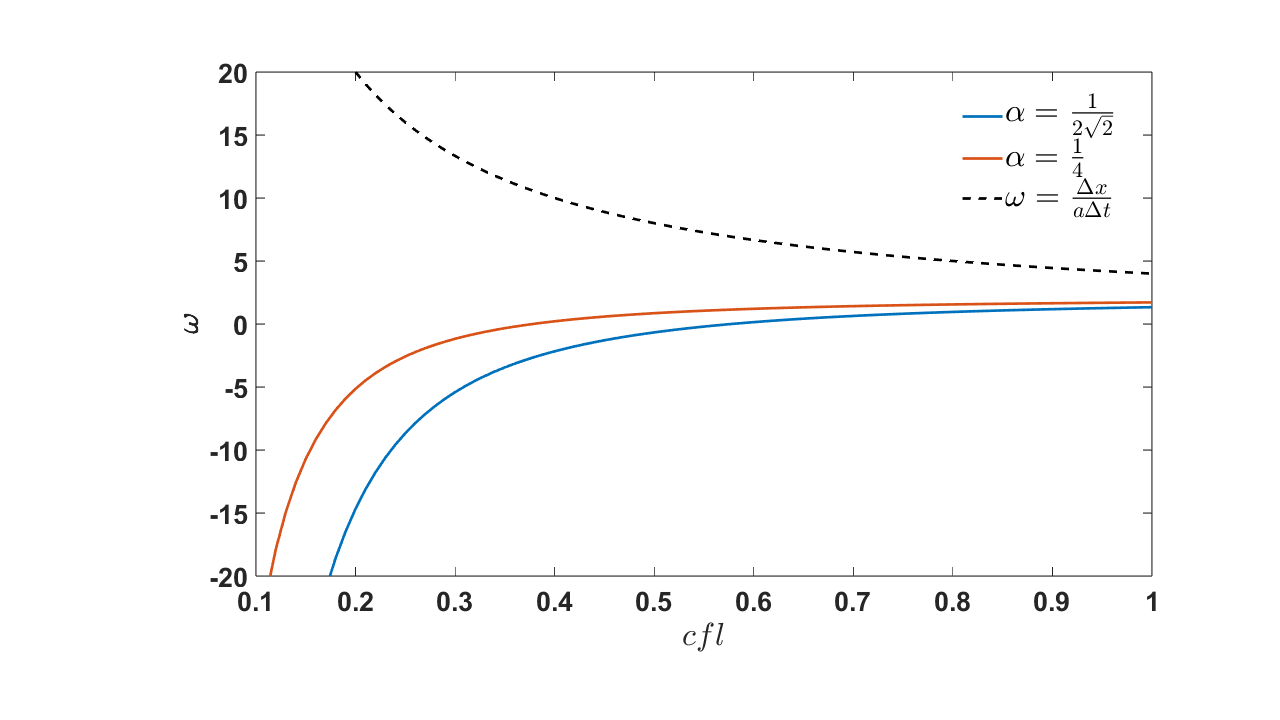} } \\
  \mbox{
          \makebox[0.47\textwidth][c]{ (a)}
  }   
    \mbox{
          \makebox[0.47\textwidth][c]{ (b)}
  }   
\label{fig:omega_vs_cfl}
\caption{(a) The flux correction factor, $\omega$, plotted versus (a) the nodal point distribution parameter, $\alpha$,  and (b), the $cfl$, to 
yield a zero diffusion coefficient in the ME in \eqref{eq:SL2_LF_general_replace_dt}. The plots in (b) are for a Chebyshev nodal distribution with $\alpha= \frac{1}{2\sqrt{2}}$ and an equidistant nodal distribution with $\alpha=\frac{1}{4}$. For reference, the trend for the upwind scheme (dashed lines), for which the diffusion coefficient is positive, is also plotted according to $\omega=\frac{\Delta x}{a \Delta t} = \frac{\Delta x}{a cfl |\Delta x|_{min}/a} = \frac{1}{cfl(1/2-\alpha)}$.}  
\label{fig:omega_vs_alpha_cflmin}
\end{figure}

% The  dependencies of the zero diffusion coefficient on $\alpha$, $\omega$, and $cfl$  are summarized   in Figures \ref{fig:omega_vs_alpha_cflmin}a) and b), respectively. While we plot $\omega$ on the vertical axis for convenience, it should not be interpreted as the variable that is dependent on $\alpha$ and $cfl$. Rather, all three variables are the independent variables to the zero diffusion coefficient and thus can be varied independently to yield the scheme with a zero first-order diffusion coefficient.   For reference, the trend for the upwind scheme, for which the diffusion coefficient is positive, is also plotted according to $\omega=\frac{\Delta x}{a \Delta t} = \frac{\Delta x}{a cfl |\Delta x|_{min}/a} = \frac{1}{cfl(1/2-\alpha})$. 
% With a decreasing value of $cfl$, the figures shows that $\omega$ reduces and can even be negative to have a zero coefficient on the first order diffusive term. This trend can be understood by inspecting (\ref{eq:SL2_zerodiffusion}); the first term is positive and increases with decreasing $cfl$, and thus the third term, $\omega$, must decrease accordingly to satisfy the equality. The difference in positive value of $\omega$ for the positive diffusion coefficient of the upwind scheme and the negative $\omega$ for the zero diffusion coefficient goes to infinity if the nodal distribution clusters around the center point when $\alpha$ tends to 1/2.
% In general, $\omega$ has to be larger for the Chebyshev distribution for the same $cfl$ to obtain zero diffusion.
The  dependencies of the zero diffusion coefficient on $\alpha$, $\omega$, and $cfl$  are summarized  in Figures (\ref{fig:omega_vs_alpha_cflmin}a) and b), respectively. While we plot $\omega$ on the vertical axis for convenience, it should not be interpreted as a variable that is dependent on $\alpha$ and $cfl$. Rather, all three variables are the independent variables used to zero the diffusion coefficient and thus can be varied independently to enforce that property.  For reference, the trend for the upwind scheme, for which the diffusion coefficient is positive, is also plotted using the relation $\omega=\frac{\Delta x}{a \Delta t} = \frac{\Delta x}{a cfl |\Delta x|_{min}/a} = \frac{1}{cfl(1/2-\alpha)}$. 
With a decreasing value of $cfl$, the graphs show that $\omega$ also decreases and can become negative in order to produce a zero coefficient on the first order diffusive term. This trend can be understood by inspecting (\ref{eq:SL2_zerodiffusion}); the first term is positive and increases with decreasing $cfl$, and thus the third term, $\omega$, must decrease accordingly to satisfy the equality. The difference in the positive value of $\omega$ for the positive diffusion coefficient of the upwind scheme and the negative $\omega$ for the zero diffusion coefficient goes to infinity if the nodal distribution clusters around the center point when $\alpha$ tends to $1/2$.
In general, $\omega$ has to be larger for the Chebyshev distribution for the same $cfl$ to obtain zero diffusion.

% %%%%%%%%%%%%%%%%%%% Figure %%%%%%%%%%%%%%%%%%%%%%%%%%%%%

\subsection{Generalized  analysis for $P>1$}
To generalize the analysis for $P>1$, we start with the monomial of degree $P$ at $t^n$ in the element $k$ as
\begin{eqnarray}
q_k^n(\xi) \approx Q^n(\xi) = C_1^{k,n} + C_2^{k,n} \xi + ... + C_{P+1}^{k,n} \xi^P=\sum_{m=0}^{P} C_{m+1}^{k,n} \xi^{m}.
\label{eq:SLN_initial}
\end{eqnarray}
The coefficients can be obtained by solving the system
\begin{align} 
\begin{bmatrix}
	C^{k,n}_1 \\
	C^{k,n}_2 \\
\vdots \\
    C^{k,n}_{P+1} 
\end{bmatrix} =
\begin{bmatrix}
      1 	&\xi_0 		&\dots 		&\xi^P_0 \\
      1 	&\xi_1 		&\dots 		&\xi^P_1\\
      \vdots         			   \\
      1 	&\xi_{P} 		&\dots 		&\xi^{P}_{P} \\
\end{bmatrix}^{-1}
\begin{bmatrix}
	Q_0^n     \\
	Q_1^n     \\
\vdots         \\
	Q_{P}^n      \\
\end{bmatrix},
\label{eq:Ckn}
\end{align}
where the matrix on the right-hand side is the inverse of the Vandermonde matrix, $\bar{\bar{V}}$. 
We will use  the definition of the vector of basis functions,
\begin{equation}
    \mathbf{v}=\left[ 1, x, x^2, ....x^P \right],
 \end{equation}
% it follows that
% \begin{equation}
%  \mathbf{Q}^n = \langle \mathbf{C}^n, \mathbf{v} \rangle   
% \end{equation}
to alternatively denote  the Vandermonde matrix as
\begin{align}
 \bar{\bar{V}} =  
 \begin{bmatrix}
     \mathbf{v_0} \\
      \mathbf{v_1}\\
      \vdots         			   \\
      \mathbf{v_{P}} \\
\end{bmatrix},
\end{align}
with $\mathbf{v}_i=[1 \ \ \xi_i \ \ \xi_i^2 \ \ ... \ \ \xi_i^{P+1}]$.
In matrix-vector notation, we write (\ref{eq:Ckn}) as 
\begin{equation}
\mathbf{C}^{k,n} = \bar{\bar{V}}^{-1}\mathbf{Q}^n,
\end{equation}
where the second superscript $k$ identifies the element number.
% To minimize the polynomial interpolation error, we use Chebyshev-Gauss quadrature nodes whose local coordinates are given by,
% \begin{eqnarray}
% \xi_m = \frac{\Delta x}{4}\left[ \cos{\left(-\frac{(m-1)\pi}{P} \right)} +1 \right] \ \ \ \ \ \  m=1,...,P+1
% \label{eq:CG_nodes}
% \end{eqnarray}
% \gbj{not right, check elsewhere}

In the next step of the semi-Lagrangian process, the polynomial is advected according to the advection velocity $a$,
\begin{eqnarray}
 \xi_m^\star & =&  \xi_m^n + \Delta t a,  \nonumber\\
 Q^\star_{m} & = & Q^n_{m}.
\end{eqnarray}
Then we project the advected solution back to the original basis. First, we find the updated 
coefficients as
\begin{equation}
     \mathbf{C}^{k,\star} = {\bar{\bar{V}}^\star}^{-1}\mathbf{Q}^{k,n},
    \label{eq:Cstar}
\end{equation}
where the modified Vandermonde matrix, $\bar{\bar{V^\star}}$, is determined using the advected coordinates with
\begin{align}
 \bar{\bar{V}}^\star =  
 \begin{bmatrix}
     \mathbf{v_0}^\star \\
      \mathbf{v_1}^\star\\
      \vdots         			   \\
      \mathbf{v_{P}}^\star \\
\end{bmatrix},
\end{align}
where $\mathbf{v}_m^\star=[1, \ \ \xi_m+a \Delta t \ \ \ \left(\xi_m + a \Delta t \right)^2, \ \  ..., \ \ \left(\xi_m + a \Delta t\right)^P]$.
Then the updated solution is
\begin{equation}
  \mathbf{Q}^\star = 
  \bar{\bar{V}} \mathbf{C}^{k,\star}   =     \bar{\bar{V}} {\bar{\bar{V}}^\star}^{-1}\mathbf{Q}^{k,n}.
  \label{eq:qstar_vec}
\end{equation}

For the ``least-squares fit" stage of the algorithm that corrects for boundary/interface solutions, we first define the error as
\begin{eqnarray}
  \epsilon(C_1^{n+1}, C_2^{n+1},...,C_{P+1}^{n+1}) &=& (Q_{i-1/2}^+ -  \sum_{m=0}^{P}C_m^{n+1} x_{i-1/2}^{m-1})^2  \nonumber \\  
   &+& (Q_{0}^\star \ \ \ \ \ - \sum_{m=0}^{P}C_{m+1}^{n+1} \xi_1^{m-1})^2  \nonumber  \\
    &+& (Q_{1}^\star \ \ \ \ \ - \sum_{m=0}^{P}C_{m+1}^{n+1} \xi_2^{m-1})^2  \nonumber  \\
    &+& ... \nonumber \\ 
        &+& (Q_{P}^\star \ \ \ - \sum_{m=0}^{P}C_{m+1}^{n+1} \xi_{P+1}^{m-1})^2  \nonumber  \\
     &+& (Q_{i+1/2}^+ - \sum_{m=0}^{P}C_{m+1}^{n+1} x_{i+1/2}^{m-1})^2.  \\
\end{eqnarray}
Setting the derivative of the error with respect to the coefficients $C_{m+1}^{n+1}$=0 
\begin{equation}
    \frac{\partial \epsilon(C_1^{n+1}, C_2^{n+1},...,C_{P+1}^{n+1})}{\partial C_{m+1}^{n+1}} =0, \ \ \ \ \ \ \ \ m=0,...,P
\end{equation}
leads to a system of $P+1$ equations that can be solved for $C_m^{n+1}$
\begin{eqnarray}
  0=&&(Q_{i-1/2}^+ -  \sum_{m=0}^{P}C_m^{n+1} x_{i-1/2}^{m})2 x_{i-2}^{r}  \nonumber \\  
   &+& (Q_{0}^\star \ \ \ \ \ - \sum_{m=0}^{P}C_{m+1}^{n+1} \xi_1^{m})2\xi_1^{r}  \nonumber  \\
    &+& (Q_{1}^\star \ \ \ \ \ - \sum_{m=0}^{P}C_{m+1}^{n+1} \xi_2^{m})2\xi_2^{r}  \nonumber  \\
    &+& ... \nonumber \\ 
        &+& (Q_{P}^\star \ \ \ - \sum_{m=0}^{P}C_{m+1}^{n+1} \xi_{P+1}^{m})2 \xi_{P+1}^{r}  \nonumber  \\
     &+& (Q_{i+1/2}^+ - \sum_{m=0}^{P}C_{m+1}^{n+1} x_{i+1/2}^{m})2x_{i+1/2}^{r} \ \ \ \ \ \  r=0,...,P.  
\end{eqnarray}
Rewritten in matrix-vector notation this leads to
\begin{eqnarray}
&&
\begin{bmatrix}
	C^{k,n+1}_1 \\
	C^{k,n+1}_2 \\
\vdots \\
    C^{k,n+1}_{P+1} 
\end{bmatrix} = \nonumber\\
&&
\begin{bmatrix}
      P+1 	&x_{i-1/2} + \xi_0 +\dots  + x_{i+1/2} 		&\dots 		&x_{i-1/2}^P + \xi_0^P +\dots  + x_{i+1/2}^P \\
      x_{i-1/2} + \xi_0 +\dots  + x_{i+1/2} &	x_{i-1/2}^2 + \xi_0^2 +\dots  + x_{i+1/2}^2	&\dots 		&x_{i-1/2}^{P+1} + \xi_0^{P+1} +\dots  + x_{i+1/2}^{P+1}\\
      \vdots         			   \\
           x_{i-1/2}^{P+1} + \xi_0^{P+1} +\dots  + x_{i+1/2}^{P+1} &	x_{i-1/2}^{P+2} + \xi_0^{P+2} +\dots  + x_{i+1/2}^{P+2}	&\dots 		&x_{i-1/2}^{2P} + \xi_0^{2P} +\dots  x_{i+1/2}^{2P} \\
\end{bmatrix}^{-1}  \nonumber \\
&&
\begin{bmatrix}
	-(Q_{i-1/2}^+ + Q_{0}^\star + Q_{1}^\star+ ... + Q_{P}^\star +   Q_{i+1/2}^+  )  \\    	-(Q_{i-1/2}^+ x_{i-1/2} + Q_{0}^\star \xi_0+ Q_{1}^\star \xi_1+ ... + Q_{P}^\star \xi_{P} +   Q_{i+1/2}^+ x_{i+1/2} )\\
\vdots         \\
	-(Q_{i-1/2}^+ x_{i-1/2}^P + Q_{0}^\star \xi_0^P+ Q_{1}^\star \xi_1^P+ ... + Q_{P}^\star \xi_{P}^P +   Q_{i+1/2}^+ x_{i+1/2}^P  )  \\
\end{bmatrix}.
\label{eq:Ckn+1}
\end{eqnarray}
We shall denote this short-hand in matrix-vector form as
\begin{equation}
    \mathbf{C}^{k,n+1} = \bar{\bar{A}}^{-1} \mathbf{b}.
\end{equation}

We next need to express the interface/boundary solutions in terms of the interpolant in the neighbouring elements. We start by determining the coefficients in the neighbouring elements
according to (\ref{eq:Cstar}), then use those coefficients to evaluate the solutions at the interface and apply the Lax-Friedrichs Flux, which leads to
\begin{eqnarray}
    Q_{i-1/2}^+ & = & \mathcal{G} \left(  \sum_{m=0}^{P} C_{m+1}^{k-1,\star} (-\Delta x/2)^{m}  ,\sum_{m=0}^{P} C_{m+1}^{k,\star} (-\Delta x/2)^{m}\right),  \nonumber \\
    Q_{i+1/2}^+ & = & \mathcal{G} \left(  \sum_{m=0}^{P} C_{m+1}^{k+1,\star} (\Delta x/2)^{m}  ,\sum_{m=0}^{P} C_{m+1}^{k,\star} (\Delta x/2)^{m}\right).
\end{eqnarray}
This completes the algorithm.  We can write the interface values in vector notation also as
\begin{eqnarray}
  Q_{i-1/2}^+ & = & \mathcal{G} \left( \langle \mathbf{C}^{k-1,\star}, \mathbf{v}_{i-1/2} \rangle ,  \langle  \mathbf{C}^{k,\star}  \mathbf{v}_{i-1/2} \rangle  \right), \nonumber \\
  &=&  \frac{a \left[ \langle \mathbf{C}^{k-1,\star}, \mathbf{v}_{i-1/2} \rangle +  \langle 
  \mathbf{C}^{k,\star}  \mathbf{v}_{i-1/2}   \rangle \right] }{2} - 
  \frac{\omega\Delta t}{2 \Delta x}
  \left( { \langle \mathbf{C}^{k-1,\star}, \mathbf{v}_{i-1/2} \rangle -  \langle 
  \mathbf{C}^{k,\star}  \mathbf{v}_{i-1/2}   \rangle} \right) \nonumber\\
  &=&  \mathbf{v}_{i-1/2} \left[\left( a/2 -\frac{\omega\Delta t}{2 \Delta x} \right) {\bar{\bar{V}}^\star}^{-1}\mathbf{Q}^{k-1,n} + \left( a/2 +\frac{\omega\Delta t}{2 \Delta x} \right)  {\bar{\bar{V}}^\star}^{-1}\mathbf{Q}^{k,n} \right] \\
  Q_{i+1/2}^+ & = & \mathcal{G} \left( \langle \mathbf{C}^{k,\star}, \mathbf{v}_{i+1/2} \rangle ,  \langle\mathbf{C}^{k+1,\star}  \mathbf{v}_{i+1/2} \rangle  \right) \\
    &=&  \mathbf{v}_{i+1/2} \left[\left( a/2 -\frac{\omega\Delta t}{2 \Delta x} \right) {\bar{\bar{V}}^\star}^{-1}\mathbf{Q}^{k,n} + \left( a/2 +\frac{\omega\Delta t}{2 \Delta x} \right)  {\bar{\bar{V}}^\star}^{-1}\mathbf{Q}^{k+1,n} \right], \\
\end{eqnarray}
where $\langle \cdot, \cdot\rangle$ denotes an inner product of two vectors.
Now using the definition of the vector that includes the boundary corrections according
to (\ref{eq:SLL_Qb}) as
\begin{equation}
    \mathbf{Q}^{b} =\left[ Q_{i-1/2}^+, \mathbf{Q}^\star, Q_{i+1/2}^+ \right],
\end{equation}
with $\mathbf{\xi}_b=[x_{i-1/2} \ \ \xi_0 \ \ \xi_1 \ \ ... x_{i+1/2} ]^T$, where
the superscript $T$ identifies a transpose,
we can express $\mathbf{b}$ as
\begin{align}
    \mathbf{b} = 
    \begin{bmatrix}
        - \langle \mathbf{Q}^{b} , \mathbf{\xi}_b^0 \rangle \\ 
        - \langle \mathbf{Q}^{b} , \mathbf{\xi}_b \rangle \\
         - \langle \mathbf{Q}^{b} ,  \mathbf{\xi}_b^2  \rangle \\
        \vdots \\
         - \langle \mathbf{Q}^{b} ,  \mathbf{\xi}_b^P \rangle \\
        \end{bmatrix}
        = -  \mathbf{Q}^{b}     
        \begin{bmatrix}
          \mathbf{\xi}_b^0  \ \ 
          \mathbf{\xi}_b  \ \
        - \mathbf{\xi}_b^2  \ \
        \cdots \ \
          \mathbf{\xi}_b^P \ \
        \end{bmatrix} = 
        -
          \begin{bmatrix}
         \mathbf{\xi}_b^0  \ \ 
          \mathbf{\xi}_b  \ \
          \mathbf{\xi}_b^2  \ \
        \cdots \ \
          \mathbf{\xi}_b^P \ \
        \end{bmatrix}^T  \left[\mathbf{Q}^{b} \right]^T,
\end{align}
where  we have used the notation $\xi_b^P$ to identify a vector whose
elements are taken to the power $P$ as
\begin{equation}
\mathbf{\xi}_b^P=\xi_{b,i}^P \ \ \  i=0,...,P.
\end{equation}
 Using this same notation, we can write the matrix $\bar{\bar{A}}$ as
\begin{align}
    \mathbf{A} = 
    \begin{bmatrix}
     \langle \mathbf{\xi}_b^0, \mathbf{\xi}_b^0 \rangle & \langle \mathbf{\xi}_b^{1/2}, \mathbf{\xi}_b^{1/2} \rangle & ... & \langle \mathbf{\xi}_b^{P/2}, \mathbf{\xi}_b^{P/2} \rangle \\
       \langle \mathbf{\xi}_b^{1/2}, \mathbf{\xi}_b^{1/2} \rangle & \langle \mathbf{\xi}_b, \mathbf{\xi}_b \rangle & ... & \langle \mathbf{\xi}_b^{(P+1)/2}, \mathbf{\xi}_b^{(P+1)/2} \rangle \\
       \vdots \\
       \langle \mathbf{\xi}_b^{P/2}, \mathbf{\xi}_b^{P/2} \rangle & \langle \mathbf{\xi}_b^{(P+1)/2}, \mathbf{\xi}_b^{(P+1)/2} \rangle & ... & \langle \mathbf{\xi}_b^P, \mathbf{\xi}_b^P \rangle \\
    \end{bmatrix}
\end{align}

Now we can express $\mathbf{Q}^{n+1}$ in terms of $\mathbf{Q}^n$ as
\begin{eqnarray}
    \mathbf{Q}^{k,n+1} &=& \bar{\bar{V}} C^{k,n+1}= \bar{\bar{V}} \bar{\bar{A}}^{-1} \vec{b}, \nonumber \\
    &=&  -\bar{\bar{V}} \bar{\bar{A}}^{-1}         
    \begin{bmatrix}
         \mathbf{\xi}_b^0  \ \ 
          \mathbf{\xi}_b  \ \
          \mathbf{\xi}_b^2  \ \
        \cdots \ \
          \mathbf{\xi}_b^P \ \
        \end{bmatrix}^T  \left[\mathbf{Q}^{b} \right]^T , \nonumber \\
   &=&   -\bar{\bar{V}} \bar{\bar{A}}^{-1}        
          \begin{bmatrix}
         \mathbf{\xi}_b^0  \ \ 
          \mathbf{\xi}_b  \ \
          \mathbf{\xi}_b^2  \ \
        \cdots \ \
          \mathbf{\xi}_b^P \ \
        \end{bmatrix}^T \left[ Q_{i-1/2}^+, \mathbf{Q}^\star, Q_{i+1/2}^+ \right]^T \\
   &=&   - \bar{\bar{V}} \bar{\bar{A}}^{-1}         
          \begin{bmatrix}
         \mathbf{\xi}_b^0  \ \ 
          \mathbf{\xi}_b  \ \
          \mathbf{\xi}_b^2  \ \
        \cdots \ \
          \mathbf{\xi}_b^P \ \
        \end{bmatrix}^T \nonumber \\
   && \left[ \mathbf{v}_{i-1/2} \left[\left( a/2 -\frac{\omega\Delta t}{2 \Delta x} \right) {\bar{\bar{V}}^\star}^{-1}\mathbf{Q}^{k-1,n} + \left( a/2 +\frac{\omega\Delta t}{2 \Delta x} \right)  {\bar{\bar{V}}^\star}^{-1}\mathbf{Q}^{k,n} \right], \right. \nonumber \\
   && \ \ \ \ \ \ \ \ \ \  \left[ \bar{\bar{V}} {\bar{\bar{V}}^\star}^{-1}\mathbf{Q}^{k,n}  \right], \nonumber \\
   &&  \ \left. \mathbf{v}_{i+1/2} \left[\left( a/2 -\frac{\omega\Delta t}{2 \Delta x} \right) {\bar{\bar{V}}^\star}^{-1}\mathbf{Q}^{k,n} + \left( a/2 +\frac{\omega\Delta t}{2 \Delta x} \right)  {\bar{\bar{V}}^\star}^{-1}\mathbf{Q}^{k+1,n} \right]  \right]^T.
\end{eqnarray}
This can be written in shortened notation as
\begin{eqnarray}
    \mathbf{Q}^{k,n+1} &=&  - \bar{\bar{V}} \bar{\bar{A}}^{-1}  M_{\xi_b} \left[M_k^T \mathbf{Q}^{k,n}   + M_{k-1}^T\mathbf{Q}^{k-1,n}  +  M_{k+1}^T\mathbf{Q}^{k+1,n}\right],
    \label{eq:MxV}
\end{eqnarray}
with the size of all matrices $M$ being $P\times(P+2)$ as
\begin{eqnarray}
M_{\xi_b} & = &  
\begin{bmatrix}
         \mathbf{\xi}_b^0  \ \ 
          \mathbf{\xi}_b  \ \
          \mathbf{\xi}_b^2  \ \
        \cdots \ \
          \mathbf{\xi}_b^P \ \
\end{bmatrix}^T,  \nonumber \\
M_{k} & = & \left[
\mathbf{v}_{i-1/2} \left[\left( a/2 +\frac{\omega\Delta t}{2 \Delta x} \right){\bar{\bar{V}}^\star}^{-1} \right]                                \ \ \  
\left[ \bar{\bar{V}} {\bar{\bar{V}}^\star}^{-1}  \right]                 \ \ \  \mathbf{v}_{i+1/2}\left[\left( a/2 -\frac{\omega\Delta t}{2 \Delta x} \right){\bar{\bar{V}}^\star}^{-1} \right] \right], \nonumber \\
M_{k-1} & = & \left[ 
\mathbf{v}_{i-1/2} \left[\left( a/2 -\frac{\omega\Delta t}{2 \Delta x}    \right){\bar{\bar{V}}^\star}^{-1} \right]                               \ \ \ 
\bar{\bar{V}}*0                                \ \ \             \ \     \ \ \ \ \  \   \ \ \ \ \ \  \ \ \ \ \  
\mathbf{v}_{i-1/2}*0 \ \ \ \ \ \ \ \ \ \ \ \ \ \right], \nonumber \\
M_{k+1} & = & \left[      \ \ \ \ \ \  \ \ 
\mathbf{v}_{i+1/2}*0  \ \ \ \ \ \  \ \ \ \ \ \ \ \ \ \ \ \ \ \ \ \
\bar{\bar{V}}*0                                                   \ \ \ \ \ \  \ \ 
\mathbf{v}_{i+1/2} \left[\left( a/2 +\frac{\omega\Delta t}{2 \Delta x}    \right){\bar{\bar{V}}^\star}^{-1} \right]                              
\right].
\label{eq:M}
\end{eqnarray}
%This matrix-vector system should lends itself to Fourier/eigenvalue analysis.
The updated solution in (\ref{eq:MxV}) through a recursive relation and a single, recursive
matrix-vector operation which lends itself to further analysis.

For the sake of a simplified analysis, it is common to take periodic
boundary conditions \cite{GassnerKopriva11}. If we assume periodic boundary conditions then 
$\mathbf{Q}^{k+1,n}=\mathbf{Q}^{k+1,n}$ and
$\mathbf{Q}^{k-1,n}=\mathbf{Q}^{k+1,n}$. In that case
$M_{k-1}$ and $M_{k+1}$ are
\begin{eqnarray}
M_{k-1} & = & \left[ 
\mathbf{v}_{i+1/2} \left[\left( a/2 -\frac{\omega\Delta t}{2 \Delta x}    \right){\bar{\bar{V}}^\star}^{-1} \right]                               \ \ \ 
\bar{\bar{V}}*0                                \ \ \             \ \     \ \ \ \ \  \   \ \ \ \ \ \  \ \ 
\mathbf{v}_{i-1/2}*0 \ \ \ \ \ \ \ \ \ \ \ \ \ \right], \nonumber \\
M_{k+1} & = & \left[      \ \ \ \ \ \  \ \ 
\mathbf{v}_{i+1/2}*0  \ \ \ \ \ \  \ \ \ \ \ \ \ \ \ \ \ \ \
\bar{\bar{V}}*0                                                   \ \ \ \ \ \  \ \ 
\mathbf{v}_{i-1/2} \left[\left( a/2 +\frac{\omega\Delta t}{2 \Delta x}    \right){\bar{\bar{V}}^\star}^{-1} \right]                              
\right],
\label{eq:Mperiodic}
\end{eqnarray}
which can be combined with $M_k$  so that
\begin{eqnarray}
M_{k} & = & \left[
\left[
\frac{a}{2}                       \left[\mathbf{v}_{i-1/2} + \mathbf{v}_{i+1/2}\right] +
\frac{\omega\Delta t}{2 \Delta x} \left[\mathbf{v}_{i-1/2} - \mathbf{v}_{i+1/2}\right] 
\right]
\left[
{\bar{\bar{V}}^\star}^{-1} 
\right]  
\right. \nonumber \\                 
&&  \ \ \left[ \bar{\bar{V}} {\bar{\bar{V}}^\star}^{-1}  \right]        \nonumber \\   
&& \left. \ \
\left[
\frac{a}{2}                       \left[\mathbf{v}_{i-1/2} + \mathbf{v}_{i+1/2}\right] +
\frac{\omega\Delta t}{2 \Delta x} \left[\mathbf{v}_{i+1/2} - \mathbf{v}_{i-1/2}\right] 
\right]
\left[
{\bar{\bar{V}}^\star}^{-1} 
\right] 
\right]. 
\end{eqnarray}
Then (\ref{eq:MxV}) simplifies to
\begin{eqnarray}
    \mathbf{Q}^{k,n+1} &=&  - \bar{\bar{V}} \bar{\bar{A}}^{-1}  M_{\xi_b} M_k^T \mathbf{Q}^{k,n}.
    \label{eq:MxVperiodic}
\end{eqnarray}
% \gbj{this may be more efficient then the Lagrange polynomial based algorithm as the matrix on  the right hand-side of
% \ref{eq:Ckn+1} can be predetermined and the vector depends on $C^\star$ and it can potential be simplified using (\ref{eq:Cstar}) as inner product preventing the need for an inverse. The determination of $C^\star$ does require an inverse, but it of a Vandermonde which can be determined a priori also. Can it be done in 2D?}
To obtain periodicity we could also just set $\mathbf{Q}^{k+1,n}=0$ and
$\mathbf{Q}^{k-1,n}=0$ in (\ref{eq:MxV}).

\subsubsection{The Modified Equation and Convergence Rates}
To  derive the Modified Equation, we must decide on the nodal point distribution. For simplicity, we will consider equidistant local node coordinates
\begin{equation}
    \xi_m = (m+1) \frac{\Delta x}{P+2} \ \ \ \ \  \ \ m=0,...,P. 
\end{equation}
These local coordinates coincide with the global coordinates of the $k^{th}$ element as
\begin{equation}
    x_{m-P/2}^k = -\frac{\Delta x}{2} + \xi_m \ \ \ \ \ \  m=0,...,P,
\end{equation}
for an even $P$ if the center nodal point is chosen at $\xi_{P/2+1}=0$.
For the sake of this analysis we determine the coefficients $C_m^{k-1,\star}$, $C_m^{k,\star}$ and $C_m^{k+1,\star}$ using the global coordinates
\begin{eqnarray}
    x_{m-P/2-P-1}^{k-1} &=& -\frac{3\Delta x}{2} + \xi_m \ \ \ \ \ \  m=0,...,P \\
    x_{m-P/2}^k &=& -\frac{\Delta x}{2} + \xi_m \ \ \ \ \ \           m=0,...,P \\
    x_{m-P/2+P+1}^{k+1} &=& \frac{\Delta x}{2} + \xi_m \ \ \ \ \ \    m=0,...,P
\end{eqnarray}
Then we use the coefficients in (\ref{eq:Ckn+1}) to determine the solution at $x=0$ as
\begin{equation}
    q_k(x=0,t^{n+1}) = Q_0 = C_1^{k,n+1} 
    \label{eq:SLP_CenterNode}
\end{equation}
Then a Taylor series analysis can be performed at this center location to obtain the modified equation. 
We have verified that the procedure just outlined leads to the MEs equations derived
for $P=0$ and $P=1$ in the previous sections. 
% \gbj{might there be point to determining the Modified Equation at other/multiple collocation points? Can they be summed to obtain modal understanding of the ME similiar to eigenvalue analysis?}

For $P=2$ we find that the ME is
\begin{align}
{\qti}+{\qxi}={\frac {3{{\it \dx}}^{2}}{70} \left( {{\it cfl}}^{2}
\omega+{\frac {35\,{{\it cfl}}^{2}}{9}}+{\frac {\omega}{16}}-{\frac{35
}{144}} \right) } \qxxxi + ...
\end{align}
The leading order term is second-order, i.e. the convergence rate is the same as $P$. We find the same order of accuracy for other nodal distributions, including the Chebyshev Gauss quadrature nodes, but the coefficients on the truncation terms are dependent on the nodal point distribution consistent with the discussion for the $P=1$ method above.
It turns out that in the derivation of MEs for other $P$, the semi-Lagrangian method is $P^{th}$ order accurate,  i.e., the leading order truncation term is of order $P$.  
To underscore, we give  the leading order term in the ME for $P=4$
    \begin{equation}
        \left(5((\omega - 231/25)cfl^4 + ((5 \omega)/12 + 77/60)cfl^2 + \omega/162 - 77/2700)\right)\frac{\dx^4}{5544}\frac{\partial^5 Q_i^n}{\partial x^5},
    \end{equation}
which is fourth order. Note that the leading order coefficient is once again
dependent on $\alpha$, $\omega$ and $cfl$ and that is possible to remove
one of the leading terms by balancing these values so that the coefficient is zero.

While the ME provides insight
into the leading order behavior of convergence, diffusion and dispersion, it was already shown
to be insufficient to understand stability, diffusion and dispersion comprehensively in the $P=0$
case. To this end, we need to understand the dispersion and diffusion for all wave content.

\subsubsection{The Diffusion and Dispersion Relations}

Following  \cite{CostaSherwinPeiro15} we use the Modified Equation to  find the approximate dispersion and diffusion relations by considering a solution of the form
\begin{equation}
   \phi = e^{(r-is)t} e^{i\kappa x} ,
   \label{eq:phi_Fourier}
\end{equation}
with $i=\sqrt{-1}$. Substituting the solution ansatz
into the governing advection equation (\ref{eq:adv_pde}) leads to
\begin{equation}
    \left\{(r-i st) + i\kappa \right\}e^{(r-is)t} e^{i\kappa x} =0 .
\end{equation}
When (\ref{eq:phi_Fourier}) is substituted
into the Modified Equation, then the high order spatial derivatives of the truncated terms change
the expression to 
\begin{equation}
    \left\{(r-i st) + i\kappa^\star \right\}e^{(r-is)t} e^{i\kappa x} =0,
\end{equation}
where $\kappa^\star$ is the effective wave number that can be directly compared to $\kappa$ of the original advection equation. Costa \textit{et al.} \cite{CostaSherwinPeiro15} report that the more expansion terms are accounted for in the Modified Equation and $\kappa^\star$, the more accurate the dispersion relation is in approximating a formal eigenvalue analysis for the entire wave number spectrum of an Eulerian grid based discontinuous Galerkin method. 

\begin{figure}
    \centering
    \includegraphics[width=0.7\textwidth]{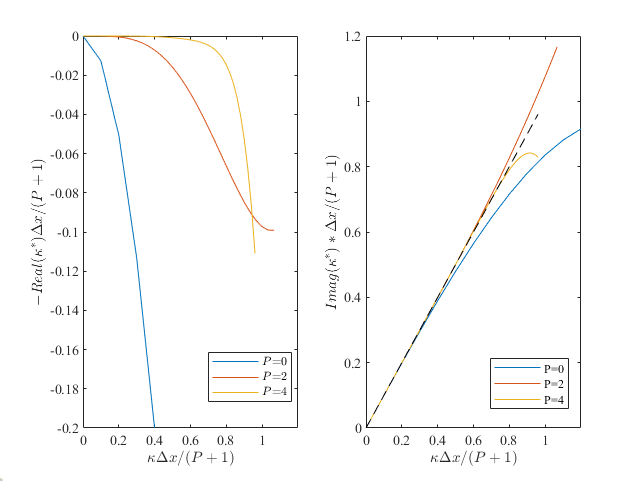}  
    \hspace{3cm}
      \mbox{
          \makebox[0.3\textwidth][c]{ (a)}
  }   
    \mbox{
          \makebox[0.3\textwidth][c]{ (b)}
  } 
    \caption{Diffusion (a) and dispersion (b) relations for $P=0,2$, and $4$.}
    \label{fig:dispersiondiffusion}
\end{figure}

Fig. \ref{fig:dispersiondiffusion} plots
the real and imaginary part of $\kappa^\star$ found using thirteen Taylor series expansion terms for $P=0$, $P=2$ and $P=4$ with an equidistant nodal distribution.
Recognizing that $-Real({\kappa^\star)}=-s/a$ 
and $Imag({\kappa^\star}=r/a$), these plots can be interpreted to understand the diffusive and dispersive nature of the semi-Lagrangian method. 

With increasing $P$ the method remains diffusive for all $\kappa$, but less so since the drop-off in the diffusion occurs at higher $P$. 
The dispersion relation for $P=0$ shows a known trend: For higher $P$, the semi-Lagrangian has a near zero dispersion as can be expected for a method that essentially traces the characteristics. For $P=2$ and $P=4$ the effective wave number is slightly higher than the exact wave number. At higher wave numbers the number of truncation terms significantly affects the accuracy of the effective dispersion relation.

\subsubsection{Stability}
The updated solution in (\ref{eq:MxVperiodic}) is recurrent and the matrix  $\bar{\bar{N}}=- \bar{\bar{V}} \bar{\bar{A}}^{-1}  M_{\xi_b} M_k^T $  is the  amplification matrix that can be analyzed by determining its eigenvalues of
% \gbj{From here a bit of help would be appreciated}.
% To determine the eigenvalues one considers the solution $\mathbf{Q}^{k,n}$  expressed as  a Fourier series
% \begin{equation}
%     \mathbf{Q}^{k,n} =  \exp^{i\left(\kappa \mathbf{x} -\omega t^n \right)}.
% \end{equation}
% Nyquist frequency, $\kappa>\kappa_{min}$ and $\omega$
\begin{equation}
    det\left[ \bar{\bar{N}} - \bar{\bar{\Lambda}} \bar{\bar{I}} \right]  = 0.
\end{equation}
The largest eigenvalue of $\bar{\bar{\Lambda}}$ represents the spectral radius and its absolute value $|\lambda_{max}|<1$, should be less than unity to ensure stability. 

The conditioning of the matrices  $\bar{\bar{V^\star}}$ and $\bar{\bar{A}}$ contribute to the magnitude of the spectral radius.
Upon inspection, $A$ appears to be diagonally dominant. But because most terms on the diagonal have terms of $\Delta x$ to some power in the denominator and because the first entry of the matrix is a constant, it follows that for small $\Delta x$, $A$ may not be positive definite.  $\bar{\bar{A}}$ is a sparse matrix because of symmetry of the collocation points. 
The conditioning of the modified Vandermonde matrix, $\bar{\bar{V^\star}}$,  depends on the distance that semi-Lagrangian particles travel away from their original seeding at the collocation nodes. The  deviation from the seeded quadrature points determines the amount of extrapolation as shown  (\ref{fig:SLschematic}) and the  related condition number of the interpolant and inverse of Vandermonde matrix. The condition number  grows with the advected distance of the collocation nodes.   

\vspace{0.5cm}
\noindent
\textit{Example for P=2}
\vspace{0.5cm}

With $P=2$ and periodic conditions according to  (\ref{eq:MxVperiodic}) the eigenvalues are
\begin{equation}
\mathbf{\Lambda} =
 \left[ \begin {array}{c} -{\frac {{\alpha}^{2}-16\,\omega}{{\alpha}^{
2}+16}}\\ \noalign{\medskip}-1\\ \noalign{\medskip}-1\end {array},
 \right] 
 \end{equation}
showing a double eigenvalue with a negative unity magnitude $-1$.  It turns out that for a single periodic element these eigenvalues are complex conjugates and hence the solution to this problem is unstable. 

\begin{figure}[!ht]
    \centering
\includegraphics[width=0.45\textwidth,trim={5cm 8cm 4cm 8cm},clip]{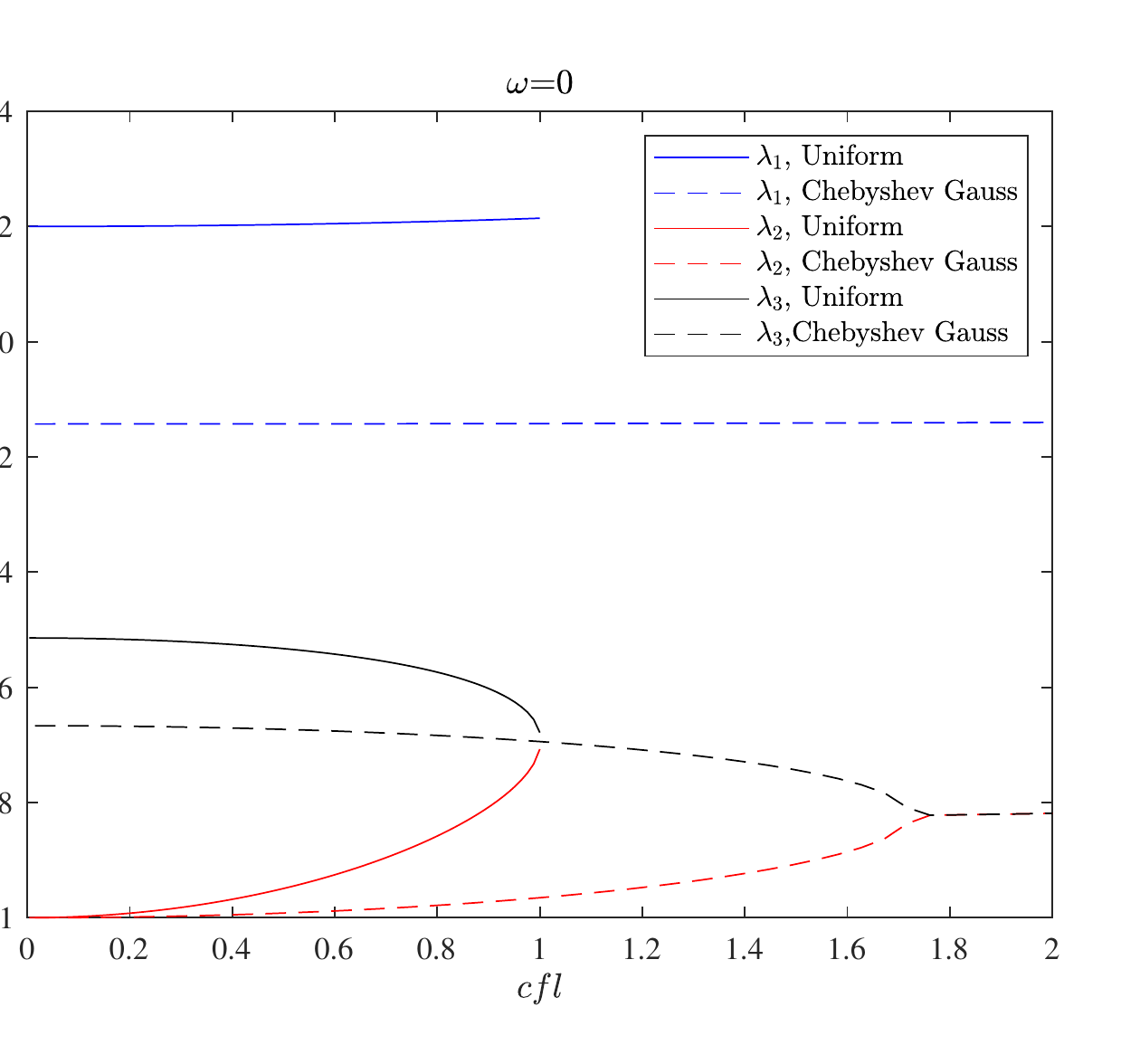}
\includegraphics[width=0.45\textwidth,trim={5cm 8cm 4cm 8cm},clip]{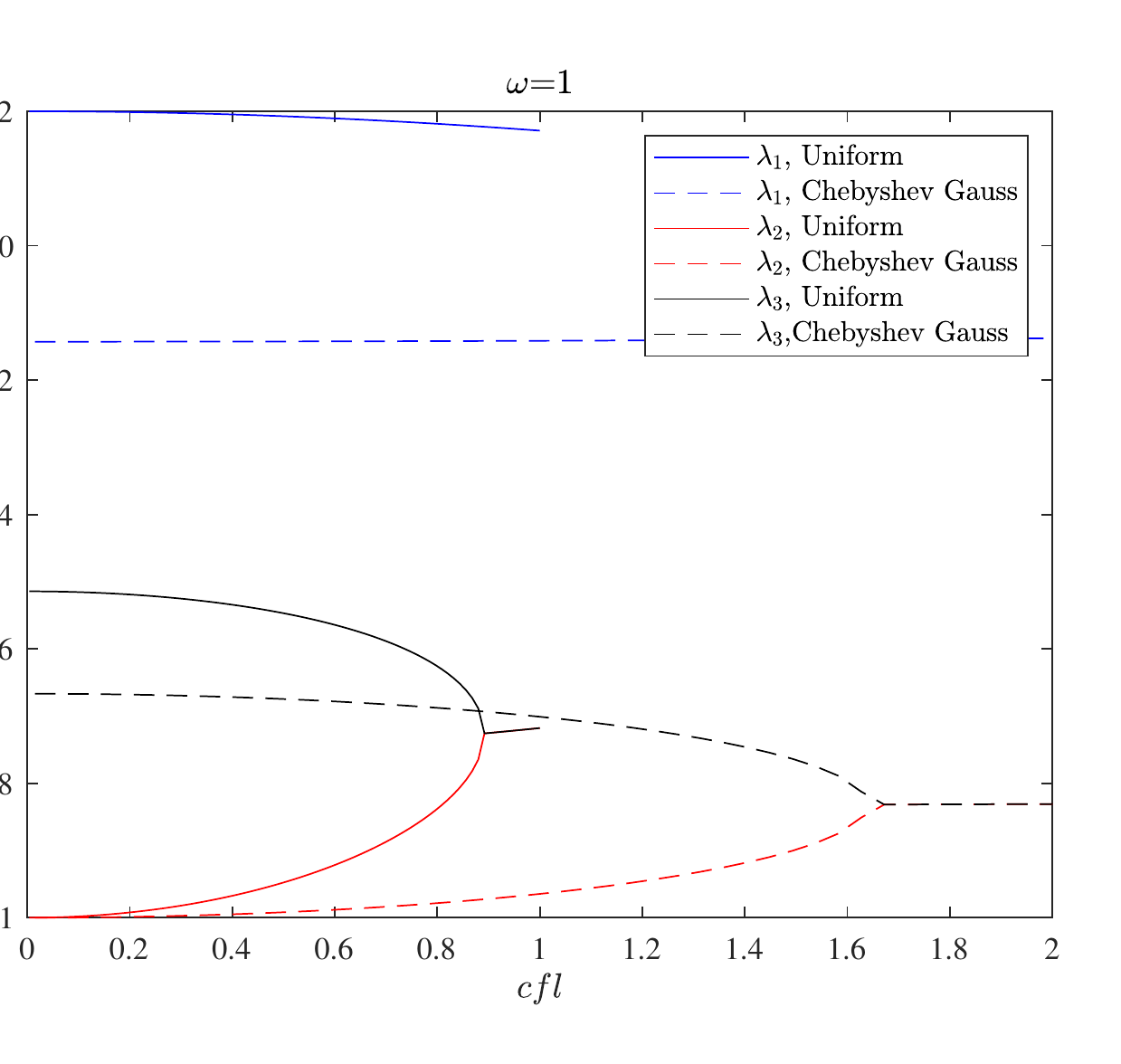}\\
\includegraphics[width=0.45\textwidth,trim={5cm 8cm 4cm 8cm},clip]{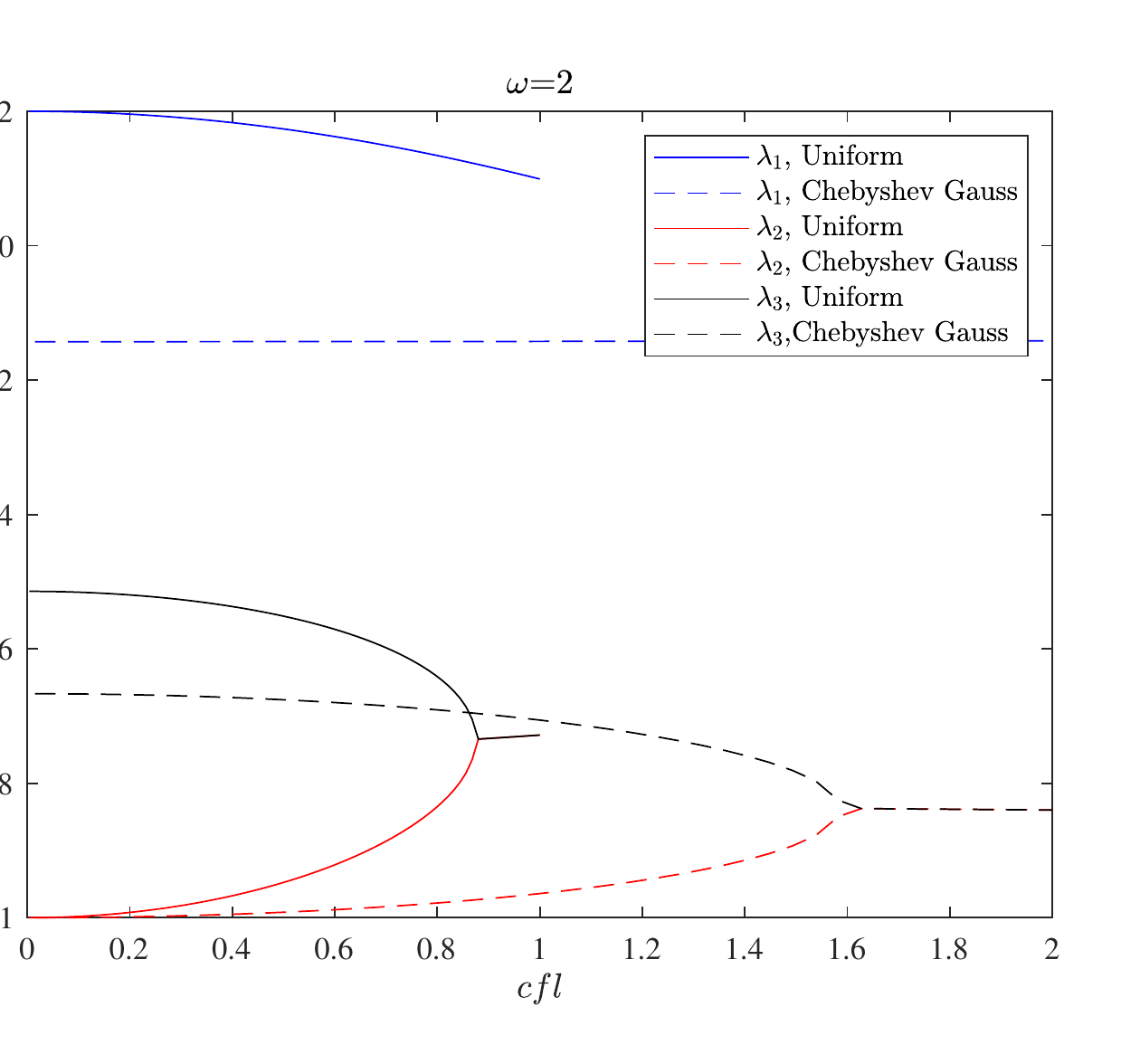}   \includegraphics[width=0.45\textwidth,trim={5cm 8cm 4cm 8cm},clip]{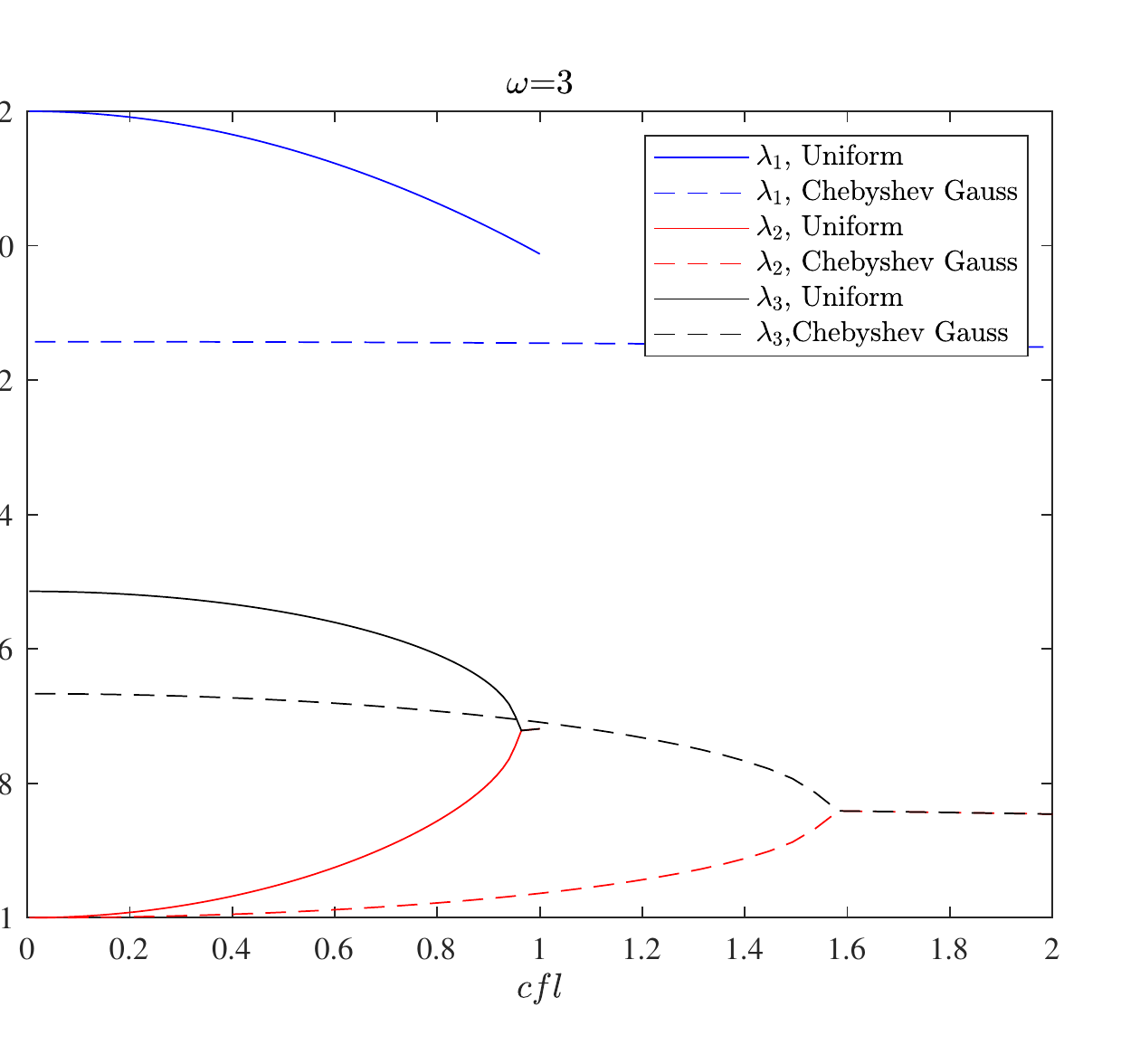}
    \caption{Eigenvalues versus $cfl$ for $\omega=0,1,2,3$ and a uniform distribution of collocation nodes and a Chebyshev quadrature node distribution with $\Delta x=0.1$.}
    \label{fig:eigenvalues}
\end{figure}

If instead of the periodic system, we impose boundary conditions by setting $\mathbf{Q}^{k+1,n}=0$ and
$\mathbf{Q}^{k-1,n}=0$ in (\ref{eq:MxV}), then we find the eigenvalues are real and unique and depend on $\cfl$ and $\omega$ as plotted in Fig. \ref{fig:eigenvalues}
for a fixed value of $\Delta x=0.1$.
Particularly, if the trends show three separate branches, then all the eigenvalues are real and their absolute
value is smaller then unity, i.e., the recursive relation and thus the semi-Lagrangian method is stable. When branches merge, the eigenvalues are complex conjugates like for the periodic system, and the solution is no longer stable. This merging of branches occurs at different $\cfl$ limits for different nodal point distributions. Chebyshev nodal distribution approximations are well-posed for $\cfl<1.6$. The $\cfl$ limit for uniform nodal increases slightly with an increased $\omega$, but is much less (smaller than unit) as compared to the limit for the Chebyshev nodal distribution. 
The eigenvalues of the lower branch tend to -1 in the limit that $\cfl$ goes to zero, suggesting
a marginally stable scheme for very small $\cfl$. This is consistent with the observation made for $P$=0 above that shows
that $\Delta x$ and $\Delta t$ should not be independently changed. 

{A von Neumann (VN) analysis can be performed  using the recurrent discrete stencil  (\ref{eq:SLP_CenterNode}) at the center node. This stencil was also used to derive the Modified Equation  for an equidistant collocation node spacing. A more precise,  stable $cfl$ criterion follows from  VN analysis as it does not require the single element periodicity assumptions that were used to find the spectral radius. The stable $cfl$ for $P=1$ with $\alpha=1/4$ according to VN is $cfl=1+\sqrt{2}$ and $cfl=\sqrt{2}$ for schemes with $\omega=\frac{\Delta x}{a(q) \Delta t}$ (upwind scheme) and $\omega=1$ (standard Lax-Friedrichs) respectively, confirming the  explicit stability  of the semi-Lagrangian method for $cfl$ greater than unity. }

\section{Numerical Tests}

Numerical tests are conducted to verify  the analysis of the semi-Lagrangian method with a Lax-Friedrichs flux correction.  
To do so,  the linear
advection equation (\ref{eq:adv_pde}) is solved with $a(q)$=constant=1 
with an  initial sine wave condition  
\begin{eqnarray}
q(x,t=0) &=& \sin(2\pi x),
\end{eqnarray}
in the domain $\left\{\Omega \in R^1: x= [0,1]\right\}$.

\subsection{Marginally resolved}

\begin{figure}[!ht]
    \centering
\includegraphics[width=0.45\textwidth]{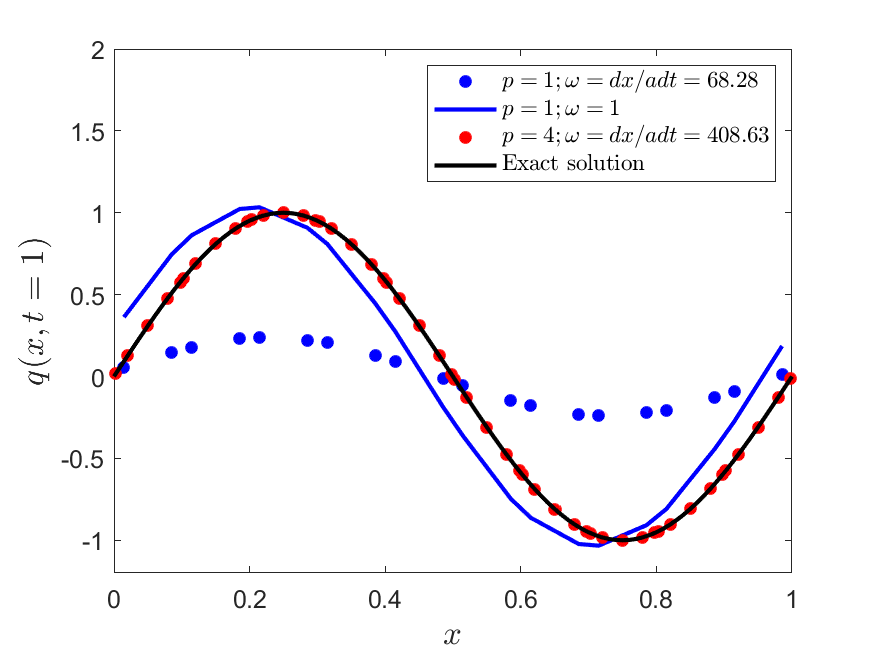}
    \caption{Numerical solution at $t=1$ for polynomial orders $P=1$ and $P=4$, flux correction factors $\omega=1$ and $\omega= dx/a dt$, ten elements, $N$=$10$, and  $cfl=\frac{1}{10}$ compared to the exact solution, $q(x,t=1)=\sin(2\pi (x-1))$.}
    \label{fig:results_marginally}
\end{figure}

We start by performing  marginally resolved simulations to illustrate  the diffusive and dispersive behavior of the first-order ($P=1$) semi-Lagrangian scheme  with  a Chebyshev spacing of the collocation nodes. To this end, we set $cfl=\frac{1}{10}$ and the number of elements, $N=10$.   We consider two cases; one for which,  $\omega= \frac{dx}{a dt}=68.28$ according to an upwind flux and for another we take $\omega=1$ according to a standard Lax-Friedrichs method.  
Figure \ref{fig:results_marginally} compares  the analytical and numerical solution so determined after one time period.  Whereas the results show that the upwind scheme is dissipative (reduced amplitude), the standard LF scheme is mostly dispersive (phase shift of the wave). This is consistent
with the Modified Equation analysis, which concluded that the leading order disspative term reduces in magnitude with reduced $\omega$. Then naturally, with a reduced diffusion, the next order dispersive  terms takes over. Also plotted in Figure \ref{fig:results_marginally} is an upwind SL solution with $P=4$, which shows that the higher order scheme is  noticeably more accurate and less dissipative.

\subsection{Convergence}
\begin{figure}[!ht]
    \centering
\includegraphics[width=0.45\textwidth]
{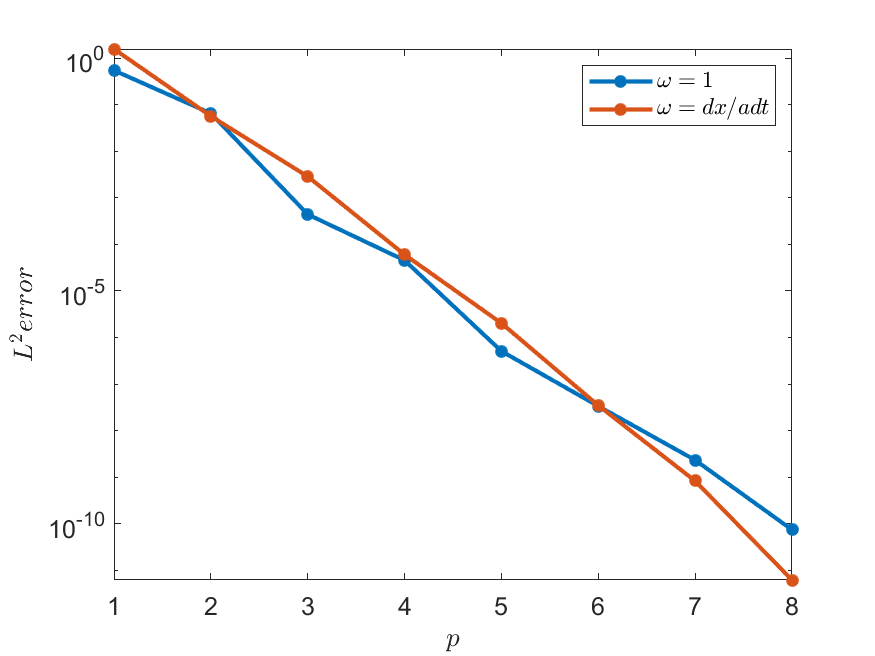}
    \caption{$L^2$ error for the semi-Lagrangian scheme plotted versus $P$ on a semi-log scale for a numerically advected sine wave over one period with $cfl$=$0.1$, $N$=$10$, and $\omega$=$ 1$ and $\omega$ = $dx/a dt$.}
    \label{fig:results_p_conv}
\end{figure}

\begin{figure}[!ht]
    \centering
\includegraphics[width=0.45\textwidth]
{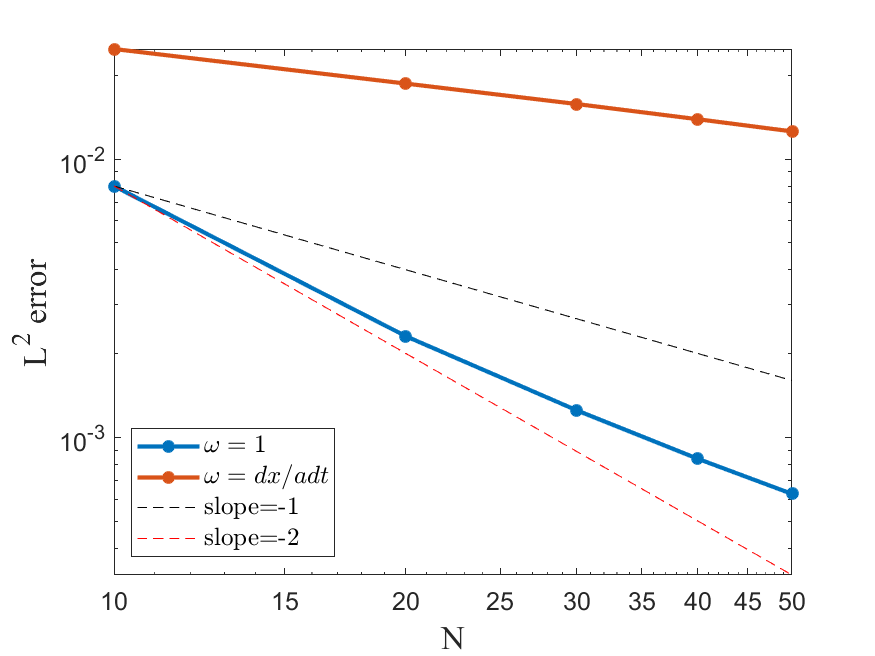}
\includegraphics[width=0.45\textwidth]
{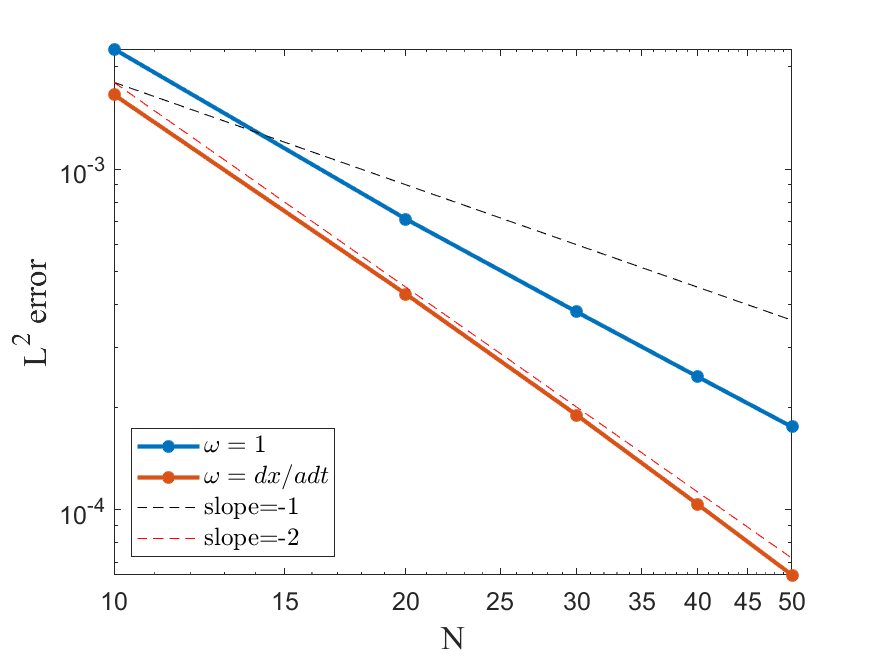}
\mbox{
\hspace{0.1cm}
\makebox[0.45\textwidth]{(a)}
 \hspace{0.05\textwidth}
\makebox[0.45\textwidth]{(b)}}
\caption{ $L^2$ error for semi-Lagrangian scheme with $P=1$ (a) and $P=2$ (b) plotted versus $N$ for several values of $\omega$ and $cfl=\frac{1}{10}$.}
    \label{fig:results_h_conv}
\end{figure}

To underscore and quantify this accuracy improvement and to validate the theoretical convergence order of $P$ as predicted by the ME analysis up to $P=8$, the convergence of the $L^2$ error versus $P$ is plotted for the uwpind SL scheme in Figure \ref{fig:results_p_conv} using  $cfl=\frac{1}{10}$ and $N=10$. The plot confirms the exponential convergence at a $P$ rate for both the upwind scheme and the scheme with $\omega=1$.  
Algebraic convergence is further confirmed to be $P^{th}$ order for $P=1$ and $P=2$  in Figure (\ref{fig:results_h_conv}). The upwind scheme exhibits a consistent
convergence rate of $P$ for a number of elements ranging from $N=10$ to $N=50$. These rates were also found for several benchmark problems in \cite{NatarajanJacobs20}.  The standard LF flux correction with $\omega=1$  appears to converge at $P+1$ for $N<20$, but  the rate reduces to $P$ for a larger number of elements.
Based on the ME analysis, this convergence rate reduction can be explained by comparing the magnitude
of the leading order diffusive term  to the next order dispersive term as follows:
for a marginally resolved numerical solution,
the leading order term in the ME is smaller for a reduced $\omega$ as compared to the upwind scheme.
So, the dispersive term may even be larger than the diffusive term at $\omega=1$.
With an increased resolution the  dispersive term asymptotically reduces faster to zero with reduced grid spacing as compared to the diffusive term  because it  is of a higher order.

\subsection{Stability}
\begin{figure}[!ht]
    \centering
\includegraphics[width=0.45\textwidth]
{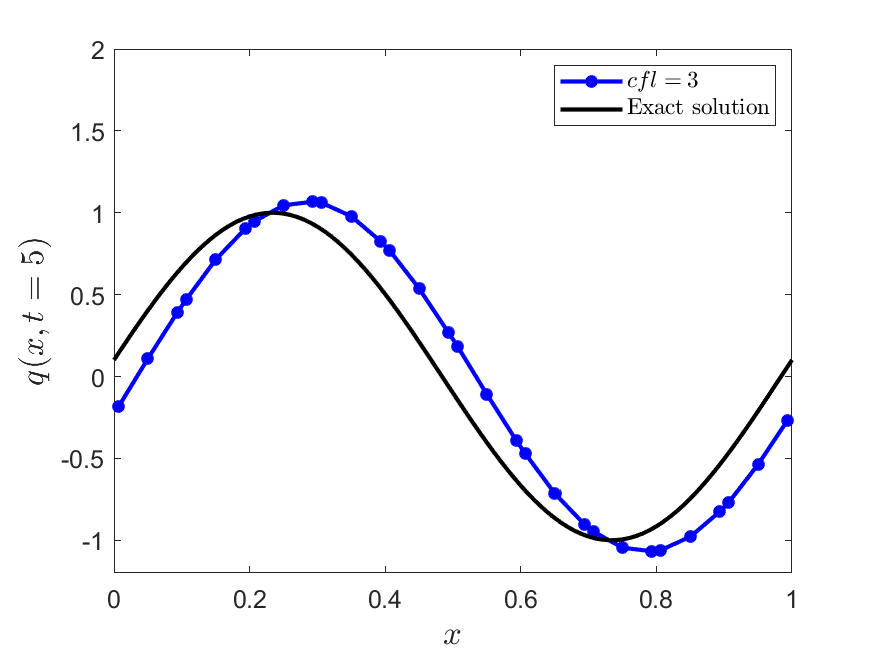}
    \caption{ The numerically advected sine solution, q(x,t=5), with $cfl$=$3$,  $P$=$2$, $N$=$10$ and $\omega$=$ dx/a dt$ as compared to the exact solution, $q(x,t=5)=\sin(2\pi (x-5))$.}
    \label{fig:results_cfl_3}
\end{figure}
To  verify the findings of the stability analysis, we simulate the sine wave with $P=2$, $N=10$ and $\omega= \frac{dx}{a(q) dt}$. The $cfl$  is increased until it is observed that the amplitude of the sine wave increases with time at which point the method is deemed unstable. Simulations are found to be  unstable  approximately for $cfl \ge  3$. Figure \ref{fig:results_cfl_3} illustrates this unstable increase of amplitude for a sine wave solution that is advected for five periods with a $cfl=3$. The expected stability limit predicted by the spectral radius analysis for a one-domain periodic case is however $cfl=1.6$.  The difference between this theoretical $cfl$ 
 and the numerically stable $cfl$ proofs that the one-domain spectral radius analysis is not sufficient to predict the exact stability limits. This observation is consistent with the conclusions made for $P=0$ and $P=1$ by comparing the results of a Von Neumann analysis for $P=0$ and $P=1$ with stability limits predicted by the ME equation analysis. Both in the ME analysis and the single domain spectral radius analysis, higher order modes or truncation terms are not accounted for that prevent the exact determination of the stable $cfl$ criterion.

\subsection{Enhanced $P$ convergence}
\begin{figure}[!ht]
    \centering
\includegraphics[width=0.45\textwidth]
{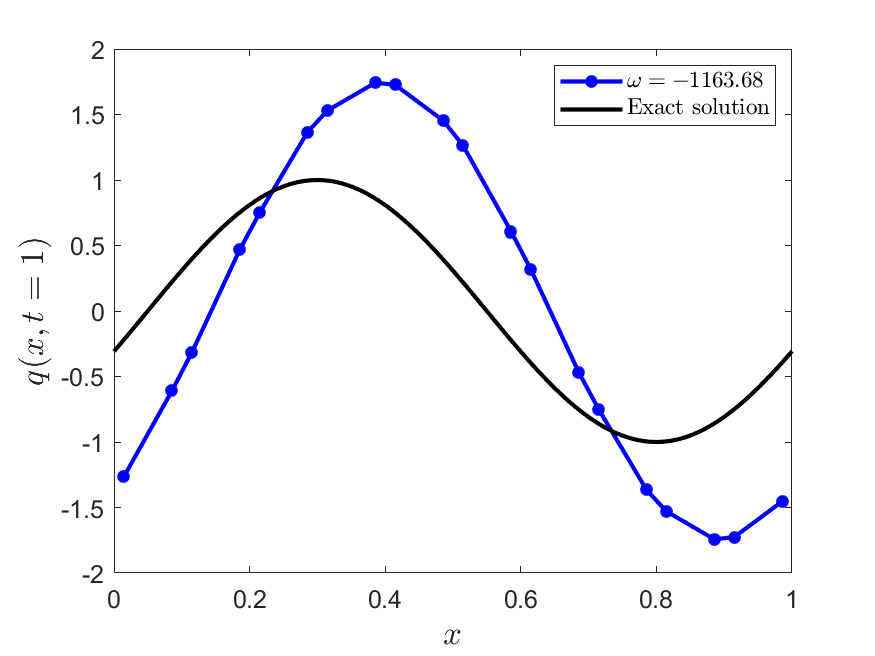}
    \caption{The numerically advected solution for $P=1$, $\omega=\omega_{\text{critical}}=-1163.68$  (yielding a zero leading order truncation term and a potential $P+1$ scheme) with $N=10$ and $cfl=\frac{1}{10}$ at $t=\frac{1}{20}$   compared to the exact solution $q(x,t=\frac{1}{20})=\sin(2\pi(x-\frac{1}{20})$.}
    \label{fig:results_a1_0}
\end{figure}
The ME equations analysis suggests that one could choose a flux correction, $\omega$,  and/or  the collocation point distribution (e.g. according to $\alpha$ for $P=1$) such that the leading order term is zero and the convergence rate increases by one order. 
% The example of the $P=1$  semi-Lagrangian scheme with $\omega=1$ discussed above suggest strongly that $P+1$ convergence may be achieved.
To test this idea, we simulate the sine wave and set $\omega=\omega_{\text{critical}}$,  the value at which the leading term in the truncation error becomes zero. For $P=1$, the truncation is according to the right hand side of \eqref{eq:SL2_LF_general_replace_dt}. 

\begin{figure}[!ht]
    \centering
\includegraphics[width=0.45\textwidth]
{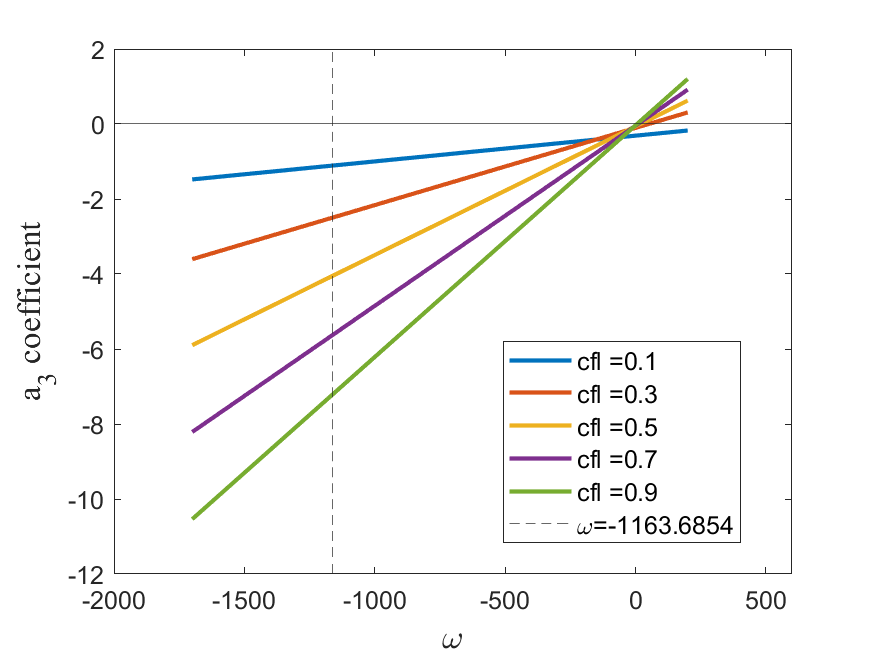}
\includegraphics[width=0.45\textwidth]
{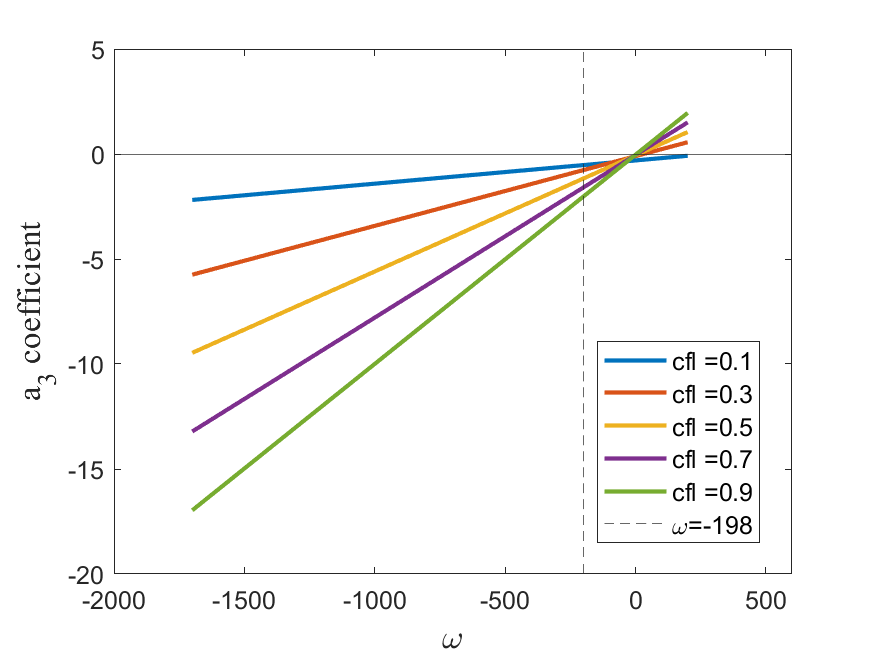}
\mbox{
\hspace{0.1cm}
\makebox[0.45\textwidth]{(a)}
 \hspace{0.05\textwidth}
\makebox[0.45\textwidth]{(b)}}
    \caption{Coefficients $a_3$ of the third-order truncation term of the Modified Equation (\ref{eq:SL2_LF_general_replace_dt})  for a $P=1$ semi-Lagrangian method for a uniform (a)
    and Chebyshev (b) collocation point distribution plotted versus $\omega$ for  several \cfl values.}
    \label{fig:results_a3_coef}
\end{figure}

Unfortunately, simulations with $cfl=0.1$ and $P=1$ using Chebyshev points with the corresponding, $\omega_{\text{critical}}=-1163.68$ are unstable as illustrated in  Figure \eqref{fig:results_a1_0}. To understand this behavior, we plot the coefficient $a_3$ of the third order truncation term
in the ME versus $\omega$ in Figure \eqref{fig:results_a3_coef}. Note that this terms is next even order diffusion term after the leading first-order diffusion term. Clearly,  the coefficient $a_3$ is  negative if $\omega = \omega_{\text{critical}}$  for a uniform and a Chebyshev nodal distribution, and tested $cfl$ values. Thus the solution is mildly anti-diffusive  and hence becomes weakly unstable.
 {We have confirmed with a Von Neumann analysis for an equidistant nodal spacing that signs in the Taylor series of the dissipation factor can alternate in the higher-order terms, and so choosing the parameters using only the leading order term need not lead to a stable scheme in the critical case.}

\section{Conclusions}

The discrete, spatio-temporal, recursive stencil is derived  for a semi-Lagrangian spectral element method that updates a piecewise
continuous interpolant according to the Lagrangian (characteristic form) of the transport equations. The  derivation
hinges on the use of the monomial form of the interpolant and the corresponding, invertible Vandermonde  matrix of
the resulting interpolant.

A Modified Equation (ME) analysis, which Taylor expands the stencil in space and time to  a single
space and time location, shows that the semi-Lagrangian
method is consistent with the PDE form of the transport equation in the limit that the grid spacing goes to zero for smooth solutions. The leading order truncation term of the ME is $P^{th}$ order, which is consistent with the numerical tests
presented in previous work \cite{NatarajanJacobs20}.

The nodal point distribution affects the coefficient of the truncation terms. It is in general possible
to increase the order of the method by one order for a linear, one-dimensional form of the equations.
It remains to be seen if this can be generalized to general element in multi dimensions.

The approximate diffusion and diffusion relations derived  by comparing the wave content of multiple truncation terms of the
ME show that the higher order
semi-Lagrangian method $P>2$ has very little dispersion at low wavenumbers, as one might expect from a semi-Lagrangian method.

An eigenvalue analysis of the recursive, algebraic form of the semi-Lagrangian method shows
that the  $\cfl$ is limited by numerical stability. The $\cfl$ limit depends
on the nodal point distributions with an element. For a Chebyshev point distribution, 
the semi-Lagrangian method is stable for values of $cfl$ significantly greater than unity.

\bibliographystyle{unsrt}

%\bibliography{thbib}
% \bibliographystyle{siamplain}

\end{document}